\documentclass[fleqn,11pt]{wlscirep}
\usepackage[misc]{ifsym}
\title{Wild Bootstrapping Rank-Based Procedures: \\ Multiple Testing in Nonparametric Split-Plot Designs}

\author[1,*]{Maria Umlauft}
\author[2]{Marius Placzek}
\author[3]{Frank Konietschke}
\author[1]{Markus Pauly}
\affil[1]{Institute of Statistics, Ulm University, Germany}
\affil[2]{Department of Medical Statistics, University Medical Center Göttingen, Germany}
\affil[3]{Department of Mathematical Science, University of Texas, Dallas, USA}
\affil[*]{Corresponding Author: Helmholtzstr. 20, 89081 Ulm, \Letter\;maria.umlauft@uni-ulm.de}
\date{\today}

\usepackage{multirow}
\usepackage{subcaption}
\usepackage{bbm}
\usepackage{bm}
\usepackage{setspace}
\usepackage{amssymb}
\usepackage{lipsum}

\newcommand{\bs}[1]{\mathbf{#1}}

\newcommand{\olw}{\overline{w}}
\newcommand{\renu}{\mathbb{R}}
 
\newcommand{\Cov}{\operatorname{{Cov}}}
\DeclareMathOperator{\Erw}{\mathbb{E}}
\DeclareMathOperator{\tr}{tr} 
\def\epsilon{\varepsilon}
\newcommand{\ep}{\epsilon}
\numberwithin{equation}{section}

\newcommand{\vwhp}{\bs{\widehat{p}}}

\newcommand{\vwhw}{\bs{\widehat{w}}}

\DeclareMathOperator{\rank}{rank}

\newcommand{\bqan}{\begin{eqnarray}}
\newcommand{\eqan}{\end{eqnarray}}
\newtheorem{definition}{{\textsc Definition}\textsc}[section]
\newtheorem{satz}[definition]{{\textsc Theorem}\textsc}
\newtheorem{corollary}[definition]{{\textsc Corollary}\textsc}

\newcommand{\bsa}{\begin{satz}}
\newcommand{\esa}{\end{satz}}
\newcommand{\bay}{\begin{array}}
\newcommand{\eay}{\end{array}}
\newcommand{\bqa}{\begin{eqnarray*}}
\newcommand{\eqa}{\end{eqnarray*}}
\def\epsilon{\varepsilon}

\keywords{Multiple Comparisons, Rank Statistics, Simultaneous Confidence Intervals, Split-Plot Designs, Wild Bootstrap}

\begin{abstract}
Split-plot or repeated measures designs are frequently used for planning experiments in the life or social sciences. Typical examples include the comparison of different treatments over time, where both factors may possess an additional factorial 
structure. For such designs, the statistical analysis usually consists of several steps. If the global null is rejected, multiple comparisons are usually performed. 
Usually, general factorial repeated measures designs are inferred by classical linear mixed models. Common underlying assumptions, such as normality or variance homogeneity are often not 
met in real data. Furthermore, to deal even with, e.g., ordinal or ordered categorical data, adequate effect sizes should be used. 
Here, multiple contrast tests and simultaneous confidence intervals for general factorial split-plot designs are developed and equipped with a novel asymptotically correct wild bootstrap 
approach.

Because the regulatory authorities typically require the calculation of confidence intervals, this work also provides simultaneous confidence intervals for single contrasts and for the ratio of different contrasts in meaningful effects.
Extensive simulations are conducted to foster the theoretical findings. Finally, two different datasets exemplify the applicability of the novel procedure.
\end{abstract}
\begin{document}
\doublespacing
%\flushbottom
\maketitle

\thispagestyle{empty}

\newpage
\singlespacing
\section{Introduction} \label{intro}

Factorial designs with repeated measures (split-plot designs) occur frequently in clinical studies or other practical applications. Usually, repeated measures experiments including more than two groups and/or different factors are 
inferred by classical linear mixed effects models postulating homoscedasticity and specific distributional assumptions (e.g. normally distributed error terms). The research question of such studies is making inference in means or other appropriate effects.
Typically, global testing procedures are applied to provide an answer to the aforementioned question. If the global null hypothesis of, e.g. \textit{no treatment effect}, \textit{no time effect} and/or \textit{no interaction effects},
in a study examining the efficacy of different treatments over a period of time, is rejected,
the typically more important questions are ``Which of the different treatments cause this rejection?'', ``Which of the different timepoints cause this rejection?''  and/or ``Which of the different interactions cause this rejection?''. Thus, multiple comparisons are performed. Finally, confidence intervals (CIs) for the corresponding effects should be calculated, since 
simultaneous confidence intervals (SCIs) for contrasts in adequate effect measures are typically required by regulatory authorities, cf. %\nocite[Ch. 5.5, p. 25]{ich9statistical}
ICH E9 Guideline (1998, ch. 5.5, p. 25)\nocite{ich9statistical}. 
%, where the confidence intervals were based on means instead of an adequate effect measure for nonparametric approaches.

Very common in practical applications are stepwise procedures using different approaches on the same data. Such procedures may lead to non-consonant test decisions, that is, the global testing procedure rejects the null hypothesis, but none of the 
individual hypotheses does or vice versa (see Gabriel, 1969\nocite{gabriel1969simultaneous}). Furthermore, the CIs and the test decision may be incompatible, since the CI may include the value of no treatment effect even if the 
corresponding null hypothesis has been rejected (Bretz et al., 2001\nocite{bretz2001numerical}). One alternative is the classical Bonferroni adjustment, which can be used to perform multiple comparisons and is also useful for computing 
compatible SCIs.  
However, it results in low power, especially if the test statistics are not independent. 
This first naive approach can often be enhanced by taking the correlation between the test statistics into account.

Parametric approaches realizing such multiple testing procedures leading to compatible SCIs were introduced in the last years. The procedures proposed by Mukerjee et al.~(1987)\nocite{mukerjee1987comparison} and 
Bretz et al.~(2001)\nocite{bretz2001numerical} are suitable in case of unpaired data assuming homogeneity and normality and also taking the correlation between the test statistics into account. The method of Bretz et al.~(2001)\nocite{bretz2001numerical} is a powerful tool 
for the computation of compatible SCIs. The clue of their procedure is the establishment of an exact joint multivariate $t$-distribution that allow the control of the familywise type-$I$ error in the strong sense 
(Hochberg \& Tamhane, 1987\nocite{hochberg1987multiple}). Extensions of Bretz et al.~(2001)\nocite{bretz2001numerical} for heteroscedastic data were given by Hasler \& Hothorn 
(2008)\nocite{hasler2008multiple} and Herberich et al.~(2010)\nocite{herberich2010robust}. Hothorn et al.~(2008)\nocite{hothorn2008simultaneous} even extend the approach of Bretz et al.~(2001)\nocite{bretz2001numerical} to general parametric models.
All of the above-mentioned publications only deal with independent data, whereas Miller (2011)\nocite{miller2011simultane} introduced a procedure for
general factorial repeated measures designs. 
A disadvantage of parametric models is that they often impose restrictive assumptions. A violation of one of the assumptions may result in a substanial loss of power and inflated type-$I$ error rates. Additionally,
ordinal, skewed, score or non-continuous data are often present in real data applications. Thus, nonparametric multiple testing methods  
and approaches to calculate compatible SCIs are needed.

One general approach for nonparametric models is to formulate the hypothesis in terms of distribution functions, see for example Akritas \& Arnold (1994)\nocite{akritas1994fully} and Akritas \& Brunner 
(1997)\nocite{akritas1997unified}.
The specific global null hypotheses were formulated as $H_0^F: \bs C  \bs F = \bs 0$ for an adequate contrast matrix 
$\mathbf{C}=(\bs c_1,\dots,\bs c_q)'$ and $\bs F=(F_{11},\dots,F_{ad})'$ and as 
$ \Omega^F :\{ \bs c_\ell' \bs F = \bs 0, \; \ell=1,\dots, q\}$ for the corresponding multiple testing problem. Drawbacks of testing hypotheses formulated in terms of distribution functions are that no easy-to-interpret treatment effects could be defined and 
that no CIs can be calculated.

It is the aim of the present work to overcome the disadvantage of no computable CIs. Thus, we propose %is to formulate the hypothesis in terms of relative effects, see Brunner \& Puri (2001)\nocite{brunner2001nonparametric}, Ryu \& Agresti (2008)\nocite{ryu2008modeling} and Brunner et al. (2016)\nocite{brunner1999rank}.
% In the case of independent observations, a pairwise effect is defined by $w_{ij}=\int F_i dF_j$, where $F_i$ is the distribution function of group $i$. Using these pairwise relative effects, nonparametric simultaneous confidence
% intervals were constructed by Munzel \& Hothorn (2001)\nocite{munzel2001unified}, Ryu (2009)\nocite{ryu2009simultaneous}, Konietschke (2009)\nocite{konietschke2009simultane} and
% Pfeiffer (2010)\nocite{pfeiffer2010simultane}. Since these pairwise relative effects are not transitive (Brown \& Hettmansperger, 2002\nocite{brown2002kruskal}),
% the application can lead to paradoxical results (cf. Efron's Dice, Thangavelu \& Brunner, 2007\nocite{thangavelu2007wilcoxon}).  
%Thus, the general idea 
to formulate the above null hypotheses in terms of $\bs p=(p_{11},\dots,p_{ad})'$, 
a vector of transitive relative treatment effects $p_{ij} = \int G dF_{ij}$, where $G= \frac{1}{ad}\sum_{i=1}^a \sum_{j=1}^d F_{ij},$ instead of $\bs F=(F_{11}, \ldots, F_{ad})'$. 
Adequate nonparametric estimates of the treatment effect are based on ranks, therefore such nonparametric approaches are often called rank-based procedures. 

In this work, we combine the approaches of Konietschke et al.~(2012)\nocite{konietschke2012rank}, where multiple contrast tests with compatible SCIs for a one-way layout with $a$ independent samples were introduced, and of 
Brunner et al.~(2018)\nocite{brunner2016rank}, where results for adequate effect measures in the univariate case were examined. In this way, we obtain rank-based multiple contrast testing procedures (MCTPs) and compatible SCIs for general factorial 
split-plot designs with adequate effect measures. Since the global testing procedure proposed in Brunner et al.~(2018)\nocite{brunner2016rank} is not asymptotically correct and the MCTP for split-plot designs may result in liberality or conservativism, a wild bootstrap approach is developed 
to circumvent these issues.
Additionally to the CIs for single contrasts, also CIs for ratios are developed in this work since such CIs are of practical interest (Dilba et al., 2004)\nocite{dilba2004simultaneous}. For example, testing for
non-inferiority of different treatment groups against a control group may be much easier if the test problem of non-inferiority margins is formulated as percentage changes.

Throughout this article, the following notation is used: The $d$-dimensional unit matrix is denoted by $\bs I_d$ and the ($d \times d$)-dimensional matrix of ones by $\bs J_d = \bs 1_d \bs 1'_d$, where $\bs 1_d = (1, \ldots, 1)'$ describes 
the $d$-dimensional column vector of ones. Furthermore, $\bs P_d = \bs I_d - (\frac{1}{d}) \bs J_d$ is the so-called $d$-dimensional centring matrix. The Kronecker product of matrices is denoted by the symbol $\otimes$.

The paper is organized as follows: In the next section, we introduce the notation and define the underlying statistical model and determine the asymptotics. In Section~\ref{sec:test}, the test statistics for the global null and the multiple contrast testing 
prodecure are introduced. Additionally, a wild bootstrap approach is described. In Section~\ref{sec:CI}, CIs for ratios are constructed and simulation results are displayed in Section~\ref{simu}. 
The novel procedures are applied to two real data examples in Section~\ref{exam}. Finally, a discussion and a conclusion are given in Section~\ref{con}.
All proofs and technical details are given in the Appendix.

\section{Statistical Model and Asymptotics}  \label{mod}

To be as general as possible, a nonparametric model with independent random vectors 
\begin{equation}\label{model}
\bs X_{ik}=(X_{i1k}, \ldots, X_{idk})', \quad i=1,\ldots, a;  \; k=1,\ldots, n_i,
\end{equation}
is studied. Here, the random variables $X_{ijk} \sim F_{ij},\; j=1,\ldots,d,$ represent $d \in \mathbb{N}$ fixed repeated measures on subject $k$ in group $i$.
For convenience, the vectors in \eqref{model} are aggregated in $\bs X =(\bs X'_{11}, \ldots, \bs X'_{an_a})'$.
Similar to Brunner et al.~(2018)\nocite{brunner2016rank},
the normalized version of the distribution function \ $F_i = \frac12 (F_i^+ + F_i^-)$ (see Ruymgaart, 1980\nocite{ruymgaart1980unified}) is used to account for 
ties in the data and for dealing with non-metric data, e.g. ordered categorical data. 
For ease of notation and computation, the \textit{relative treatment effect} of distribution function $F_{ij}$ (see Brunner et al., 2018\nocite{brunner2016rank}) with respect to the unweighted mean distribution function 
$G=\frac{1}{ad} \sum_{i=1}^a \sum_{j=1}^d F_{ij}$, i.e.
\begin{eqnarray} \label{releff}
p_{ij} \ = \ \int G dF_{ij} \ = \ \olw_{\cdot \cdot ij}, \;i=1,\ldots,a \;\text{ and }\;  j=1,\ldots,d,
\end{eqnarray}
is defined via the so-called \textit{pairwise relative effects}
\begin{equation}\label{wli}
%\begin{split}
w_{rsij}= P(X_{rs 2} < X_{ij1}) + \frac12 P(X_{rs 2} = X_{ij1}) 
 = \int F_{rs} d F_{ij}, 
% \end{split}
\end{equation}
for $r, i=1,\ldots, a$ and $s,j=1,\ldots, d$. The relative treatment effect $p_{ij}$ describes the effect of group $i$ and repeated measure $j$ with respect to a randomly chosen group and repeated measures combination. %These parameters are appropriate for metric as well as ordinal data and 
In particular, it has a nice interpretation as the probability $P(X_{ij1} \le Z) = p_{ij}$, where $Z\sim G$ is independent of $X_{ij1}$, in the case of continuous data. 
Since it does not depend on sample sizes, it is a model constant which can be used to formulate adequate hypotheses and CIs, which are given in the next section.

\subsection{Hypotheses and confidence intervals}
Setting $\bs p = (p_{11}, \ldots, p_{ad})'$ and denoting an arbitrary contrast matrix of interest by $\bs C=(\bs c_1,\dots,\bs c_q)'\in \renu^{q\times ad}$,
we are interested in developing asymptotically valid tests and SCIs
for the family of hypotheses 
\begin{eqnarray}\label{H0p_mult}
 \Omega^p :\{ \bs c_\ell' \bs p = \bs 0, \; \ell=1,\dots,q\},
\end{eqnarray}
where $\bs c_\ell=(c_{\ell11},\ldots, c_{\ell ad})' \in \renu^{ad}$.
In addition, we also derive SCIs for ratios
\begin{eqnarray}\label{ratios}
 \theta_\ell = \bs c_\ell' \bs p / \bs d_\ell' \bs p, \; \ell=1,\dots,q
\end{eqnarray}
of different contrasts $\bs c_\ell, \bs d_\ell \in \renu^{ad}$.
Note, that SCIs based on relative effects ratios for \eqref{ratios} have only been studied by Munzel (2009)\nocite{munzel2009nonparametric} for non-inferiority analyses in specific three-arm trails. CIs for mean ratios were introduced 
in, e.g. Dilba et al.~(2004, 2006)\nocite{dilba2004simultaneous}\nocite{dilba2006simultaneous} and Hasler (2009)\nocite{haslerextensions}.
CIs for the global hypothesis $H_0^p: \{\bs C \bs p = \bs 0\}$ have recently been studied in Brunner et al.~(2018)\nocite{brunner2016rank} and SCIs for \eqref{H0p_mult} in Konietschke et al.~(2012)\nocite{konietschke2012rank} for the case of $d=1$. 
How to estimate the quantities introduced so far is the topic of the next section.

\subsection{Estimation}
%Now, the focus lies on the estimation of the presented effect measures. All the 
The above  quantities can be estimated by replacing the distribution function by its empirical counterparts 
\[\widehat{F}_{ij}(x)~=~\frac{1}{n_{i}}\sum_{k=1}^{n_i}c(x-X_{ijk}),\; i~=~1,\ldots,a\;\text{ and }\; j~=~1,\ldots,d.\]
Here, $c(u)$ denotes the normalized version of the counting function, i.e. $c(u)=0,\frac{1}{2},1$, when $u$ is respectively less than, equal to, or greater than 0. Plugging the empirical distribution function $\widehat{F}_{ij}$ into \eqref{wli}, 
we obtain estimators of the pairwise effects by
\[\widehat{w}_{rsij}=\int \widehat{F}_{rs}d\widehat{F}_{ij}=\frac{1}{n_r}\left(\overline{R}_{ij\cdot}^{(ij+rs)}-\frac{n_i+1}{2}\right).\]
Here, $R_{ijk}^{(ij)}$ denotes the midrank of observation $X_{ijk}$ among all $n_i$ observations of combination $(i,j)$ and, thus, $R_{ijk}^{(ij+rs)}$ is the midrank of observation $X_{ijk}$ among all $n_i+n_r$ observations of combinations $(i,j)$ and $(r,s)$.
The overlined quantities denote the averages over the dotted index. %Also in case of ties, the estimates $\widehat{w}_{rsij}$ are unbiased and $L_2$-consistent estimates of $w_{rsij}$, see Brunner et al. (2016)\nocite{brunner2016rank}. 
Finally, the estimate of the relative effect \eqref{releff} is denoted by
\begin{equation}\label{equ:hatreleff}
 \widehat{p}_{ij}=\int\widehat{G}d\widehat{F}_{ij}=\frac{1}{ad}\sum_{r=1}^a\sum_{s=1}^d\widehat{w}_{rsij}=\frac{1}{ad}\sum_{r=1}^a\sum_{s=1}^d\frac{1}{n_r}\left(\overline{R}_{ij\cdot}^{(ij+rs)}-\frac{n_i+1}{2}\right).
\end{equation}
To derive (or estimate) the whole relative effects vector $\bs p$ (or their estimators $\widehat{\bs p}$) from the pairwise effects  the following representation given in Brunner et al.~(2018)\nocite{brunner2016rank} is used \begin{align*}
                                                                                                                                                     \bs p=\bs E_{ad}\cdot\bs w \; \text{ and } \;
                                                                                                                                                     \widehat{\bs p}=\bs E_{ad}\cdot\widehat{\bs w},                                                                                                                                           
                                                                                                                                    \end{align*}
where $\bs E_{ad}=\bs I_{ad} \otimes\left(\frac{1}{ad}\cdot\bs 1'_{ad}\right)$ and $\bs w=\left(\bs w'_{11}, \ldots, \bs w'_{ad}\right)'$ with entries $\bs w_{ij}=(w_{11ij}, \ldots, w_{adij})'=\int \bs F dF_{ij}$ for $\bs F=(F_{11}, \ldots, F_{ad})'$.

\subsection{Asymptotics and the covariance matrix}
The asymptotic covariance matrix $\bs V_N$ of $\sqrt{N}\left(\widehat{\bs p}-\bs p\right)$ can be represented as $\bs V_N=\bs E_{ad}\cdot\bs S \cdot \bs E_{ad}$, where $\bs S$ denotes the asymptotic covariance matrix of $\sqrt{N}\left(\widehat{\bs w}-\bs w\right)$.
Another representation of the components of $\sqrt{N}\left(\widehat{\bs w}-\bs w\right)$ occurs from the projection method (see e.g. Brunner \& Munzel, 2000\nocite{brunner2000nonparametric}):%, for a proof in the two-sample case):
   \begin{align*}
    \sqrt{N}(\hat{w}_{rsij}-w_{rsij}) & \doteqdot \sqrt{N}\left(\frac{1}{n_i}\sum\limits_{k=1}^{n_i}\left[\widehat{F}_{rs}(X_{ijk})-\widehat{w}_{rsij}\right]-\frac{1}{n_r}\sum\limits_{k=1}^{n_r}\left[\widehat{F}_{ij}(X_{rsk})-\widehat{w}_{ijrs}\right]\right)
     =: \sqrt{N}Z_{rsij}, 
   \end{align*}
where $\doteqdot$ denotes asymptotic equivalence ($N \to \infty$) of two sequences of random variables. %on the right- and left-hand side. 
Using this representation, Brunner et al.~(2018)\nocite{brunner2016rank} show that  $\sqrt{N}\left(\widehat{\bs p}-\bs p\right)$ is asymptotically
multivariate normally distributed with expectation zero and asymptotic covariance matrix $\bs V_N=\bs E_{ad}\cdot\Cov(\sqrt{N}\bs Z)\cdot \bs E_{ad}$, where $\textbf{Z}=(\textbf{Z}_{11}', \ldots, \textbf{Z}_{ad}')'$ with 
$\textbf{Z}_{ij}=(Z_{11ij}, \ldots, Z_{adij})'$. Let $ \bm{\Sigma}= \left( \Sigma_{rs,ij}\right)_{r,s,i,j=1}^{a,d,a,d}$ be a shorter notation for $\Cov(\sqrt{N}\bs Z)$ with block-wise entries
\begin{align*}
  \Sigma_{rs,rs}&=\Cov(\sqrt{N}\bs Z_{rs})=\left(\sigma_{rs}(p,q,p',q')\right)_{p,q,p',q'=1}^{a,d,a,d}, \\
  \Sigma_{rs,ij}&=\Cov(\sqrt{N}\bs Z_{rs}, \sqrt{N}\bs Z_{ij})=\left(\sigma_{rs, ij}(p,q,p',q')\right)_{p,q,p',q'=1}^{a,d,a,d},
\end{align*}
and (co-)variances
\begin{align*}
 \sigma_{rs}(p,q,p',q') &= N\Cov(Z_{pqrs}, Z_{p'q'rs}),  \; (r,s) = (i,j), \\
 \sigma_{rs,ij}(p,q,p',q') &= N\Cov(Z_{pqrs}, Z_{p'q'ij}), \; (r,s) \neq (i,j).
\end{align*}
The explicit formulas for the variances $\sigma_{rs}(p,q,p',q')$ and the covariances $\sigma_{rs,ij}(p,q,p',q')$ are rather cumbersome and given in Appendix \ref{sec:tecdet_cov}. Depending on $r,s,i,j,p,q,p',q'$, they are linear combinations of the following quantities
 \begin{equation}\label{equ:tau}\tau_r^{(s,j)}(p,q,p',q')=\frac{1}{n_r}\Erw\left[\left(F_{pq}(X_{rs1})-w_{pqrs}\right)\left(F_{p'q'}(X_{rj1})-w_{p'q'rj}\right)\right].\end{equation}
Plugging in the empirical distribution functions and
the estimators for the pairwise effects into \eqref{equ:tau} leads to consistent estimators of $\tau_r^{(s,j)}(p,q,p',q')$ given by
$\widehat{\tau}_r^{(s,j)}(p,q,p',q')=\frac{1}{n_r(n_r-1)}\sum_{k=1}^{n_r}D_{rsk}(p,q)\cdot D_{rjk}(p',q')$,
where $D_{rsk}(p,q):=\widehat{F}_{pq}(X_{rsk})-\widehat{w}_{pqrs}$ (see Brunner et al., 2018\nocite{brunner2016rank}).
Using the quantities discussed in the previous section, the different testing procedures based on these findings are introduced next.

%\subsection{Test statistic for the global null hypothesis}

%\subsection{Bootstrap as finite sample correction}\label{sec:boot}

\section{Test Statistics}\label{sec:test}
In this section, different methods for testing the hypotheses of interest (Section~\ref{mod}) will be introduced. All of the proposed methods are based on the point estimator $\widehat{\bs p}$. 
First, the multiple 
contrast testing procedures will be defined. Thereafter, a purely global testing procedure will be described and finally, a wild bootstrap approach for both methods will be introduced.

\subsection{Multiple Contrast Testing Procedures}\label{subsec:MCTP}
Regarding the MCTP,
the idea of Konietschke et al.~(2012)\nocite{konietschke2012rank} is generalized to split-plot designs by utilizing techniques of Placzek (2013)\nocite{placzeknichtparametrische} and Brunner et al.~(2018)\nocite{brunner2016rank}. 
If $N \to \infty$ such that $N/n_i\rightarrow \kappa_i\in(0,\infty)$ and if $\bs V_N \to \bs V$ such that $\rank(\bs V_n)=\rank(\bs V)$, a test statistic for each individual hypothesis $H_{0,\ell}^p: \bs{c}_\ell'\bs{p}=0, \; \ell=1,\ldots,q$, 
of $ \Omega^p$ in \eqref{H0p_mult} is defined as
\begin{equation}T_\ell^p=\frac{\sqrt{N}\bs{c}_\ell'\left(\widehat{\bs{p}}-\bs{p}\right)}{\sqrt{\widehat{v}_{\ell\ell}}}, \; \ell=1,\ldots,q,\end{equation}
where $\widehat{v}_{\ell\ell}=\bs{c}_\ell'\widehat{\bs{V}}_N\bs{c}_\ell$ is a consistent estimator of $v_{\ell\ell}=\bs{c}_\ell'\bs{V}\bs{c}_\ell$. Since $\sqrt{N}\bs{c}_\ell'\left(\widehat{\bs{p}}-\bs{p}\right)$ is asymptotically normally distributed with mean zero and variance $v_{\ell\ell}$, it follows from an application of 
Slutsky's theorem that $T_\ell^p$ is
standard normally distributed. To construct MCTPs and SCIs, the individual test statistics $T_\ell^p$ are collected in the vector
\begin{equation}\bs T_N=\left(T_1^p, \ldots, T_q^p\right)'.\end{equation}
Again using the asymptotic normality of $\sqrt{N}\bs{C}\left(\widehat{\bs{p}}-\bs{p}\right)$ and Slutsky's theorem, one can show that $\bs T_N$ is asymptotically multivariate normally distributed with
mean zero and correlation matrix $\bs R$, where the components of $\bs R=(r_{\ell m})_{\ell,m=1}^q$ are given by $(r_{\ell m})_{\ell,m}= v_{\ell m}/\sqrt{v_{\ell\ell}v_{mm}}$ with 
$v_{\ell m}=\bs{c}_\ell'{\bs{V}}\bs{c}_m$.

With the knowledge of the asymptotic distribution of $\bs T_N$, multiple contrast tests and SCIs can be constructed: Transferring the results of Konietschke et al.~(2012)\nocite{konietschke2012rank} for the special case of $d=1$ to the present set-up, 
note that
$ \Omega^p$ and $\bs T_N$  asymptotically generate a joint testing family.
Thus, a simultaneous testing procedure (STP) can be obtained with the help of two-sided, equicoordinate $(1-\alpha)$-quantiles $z_{1-\alpha,2,\bs R}$ of ${N}(\bs 0, \bs R)$ (Bretz et al., 2001\nocite{bretz2001numerical}) given by
\[P\left(\bigcap_{\ell=1}^q\{-z_{1-\alpha,2,\bs R}\leq T_\ell^p \leq z_{1-\alpha,2,\bs R}\}\right)=1-\alpha.\]

By replacing $\bs V$ with a consistent estimator $\widehat{\bs V}_N$ (Brunner et al., 2018), estimates of $v_{\ell\ell}$ and $v_{\ell m}$ are obtained, which can be used to construct a consistent estimator $\widehat{\bs R}=(\widehat{r}_{\ell m})$ of the correlation matrix $\bs R$, 
where $\widehat{r}_{\ell m}=\widehat{v}_{\ell m}/\sqrt{\widehat{v}_{\ell\ell}\widehat{v}_{mm}}$.
Thus, the individual test hypothesis $H_{0,\ell}^p: \bs c'_\ell \bs p=  0$ is rejected if $|T_\ell^p|\geq z_{1-\alpha,2,\widehat{\bs R}}$ and asymptotic $(1-\alpha)$-SCIs for $\bs c'_\ell \bs p$ are given by
\begin{equation}
\left[\bs{c}'_\ell \widehat{\bs{p}}-z_{1-\alpha,2,\widehat{\bs{R}}}\sqrt{\frac{\widehat{v}_{\ell\ell}}{N}};\; \bs{c}'_\ell \widehat{\bs{p}}+z_{1-\alpha,2,\widehat{\bs{R}}}\sqrt{\frac{\widehat{v}_{\ell\ell}}{N}}\right], \; \ell=1,\ldots,q.
\end{equation}
Moreover, a test procedure for the global null hypothesis $H_0^p:\bs{Cp}=\bs 0=\cap_{\ell=1}^q H_{0,\ell}^p$ is given by %$z_{1-\alpha,2,\widehat{\bs R}}$ of ${N}(\bs 0, \widehat{\bs R})$: 
\begin{equation}\label{equ:globalasym}\max\{|T_1^p|, \ldots, |T_q^p|\}\geq z_{1-\alpha, 2,\widehat{\bs R}}.\end{equation}
In the next section, another testing procedure for the global null based on the work of Brunner et al.~(2018) is presented.                                                                     

\subsection{Global Testing Procedure}
Beneath~\eqref{equ:globalasym}, a different approximate test procedure for the global null hypothesis $H_0^p:\bs{Cp}=\bs 0$ has been discussed in
Brunner et al.~(2018)\nocite{brunner2016rank}. They propose the application of an ANOVA-type statistic (ATS)
\begin{equation}\label{equ:QN}
 Q_N(\bs{C})=F_N(\bs{M})=\frac{N}{\tr{(\bs{M}\widehat{\bs{V}}_N)}} \widehat{\bs{p}}'\bs{M}\widehat{\bs{p}},
\end{equation}
where $\bs M=\bs C'(\bs{CC}')^+\bs C$ is the projection matrix on the column space of $\bs C$ and $\widehat{\bs V}_N$ is a consistent estimator of $\bs V_N$. 
Since the asymptotic distribution of $Q_N(\bs M)$ under the null is non-pivotal and rather complex, the authors proposed a Box (1954)\nocite{box1954some}-type approximation $Q_N(\bs{M}) \approx \frac{\chi^2_f}{f},$
where $f$ can be estimated by $\hat{f}= \frac{[\tr(\bs{M}\widehat{\bs{V}}_N)]^2}{\tr(\bs{M}\widehat{\bs{V}}_N\bs{M}\widehat{\bs{V}}_N)}$.
This approximation, however, leads to testing procedures, which asymptotically do not have the proposed significance level $\alpha$. To solve this issue, a wild bootstrap approach leading to asymptotically correct global inference procedure
based on \eqref{equ:QN} as well as alternatives to the SCIs and MCTPs proposed in Section~\ref{subsec:MCTP} is introduced in the next step.

\subsection{The Wild Bootstrap}
\subsubsection{Global Testing Procedure}\label{sec:boot}
For $i=1,\ldots, a$ and $k=1,\ldots,n_i$, let $\ep_{ik}$ be independent and identically distributed Rademacher variables with distribution $P(\ep_{ik} = \pm 1)=1/2$. 
Using the Rademacher variables as multipliers, a wild bootstrap version of $\sqrt{N}(\vwhw - \bs w$) can be defined for $r,i=1,\ldots,a$ and $s,j=1,\ldots,d$ as %(Ordnung wie in Brunner et al.)
\begin{align}\label{equ:boot}
\begin{split}
    \sqrt{N}\widehat{\bs{w}}^\varepsilon &=\sqrt{N}\left(\widehat{w}_{rsij}^\varepsilon\right)_{r,s,i,j}\\ 
    &=\sqrt{N}\left[\frac{1}{n_i}\sum_{k=1}^{n_i} \varepsilon_{ik}\left(\widehat{F}_{rs}(X_{ijk})-\widehat{w}_{rsij}\right)-\frac{1}{n_r}\sum_{k=1}^{n_r} \varepsilon_{rk}\left(\widehat{F}_{ij}(X_{rsk})-\widehat{w}_{ijrs}\right)\right].
\end{split}
\end{align}
Note, that we utilize identical Rademacher variables for each repeated measure ($j=1,\ldots,d$) to mimic the correct covariance structure in the limit, see Theorem \ref{asy1:rm} below.
Utilizing the expression $\widehat{\bs p}=\bs E_{ad} \cdot \widehat{\bs w}$, a wild bootstrap version of $\sqrt{N}(\vwhp - \bs p)$ is obtained by $ \vwhp^\ep = \bs E_{ad} \cdot \vwhw^\ep$.
%To receive an adequate testing procedure it must be shown that $\sqrt{N}(\vwhp - \bs p)$ and $\sqrt{N} \vwhp^\ep$ have the same distribution. Thus, the 
The following theorem ensures that the distribution of $\sqrt{N} \vwhp^\ep$ always approximates the null
distribution of $\sqrt{N}(\vwhp - \bs p)$.

\bsa\label{asy1:rm}
If $N/n_i\rightarrow \kappa_i\in(0,\infty)$, %the distribution of 
the random vector $\sqrt{N} \vwhp^\ep$, conditioned on the data, has asymptotically (as $N \to \infty$) a multivariate normal distribution with mean zero and covariance matrix $\bs V_N= \bs E_{ad} \cdot \Cov(\sqrt{N}\bs Z) \cdot \bs E_{ad}$ 
in probability, i.e. coincides with the asymptotic distribution of $\sqrt{N}(\widehat{\bs p}-\bs p)$.
\esa

As first application, we obtain a wild bootstrap test in the statistic of $Q_N(\bs{M})$ for the global null $H_0^p:\bs C \bs p=\bs 0$. To calculate adequate critical values, we define a wild bootstrap version of the ANOVA-type statistic \eqref{equ:QN} as
$Q_N^\epsilon(\bs M)= \frac{N}{\tr(\bs M \widehat{\bs V}_N^\epsilon)}\widehat{\bs p}^{\epsilon'} \bs M \widehat{\bs p}^\epsilon$,
where $\widehat{\bs V}_N^\epsilon$ denotes the covariance matrix based on the wild bootstrap samples \eqref{equ:boot}. 

\begin{corollary}\label{cor:ATS}
If $N/n_i\rightarrow \kappa_i\in(0,\infty)$, the distribution of the wild bootstrap version of the ANOVA-type test statistic $Q_N^\epsilon(\bs M)$, conditioned on the data, 
 always approximates the null distribution of $Q_N(\bs M)$ as $N \to \infty$ in probability, i.e. for every $\bs p \in \left[0,1 \right]^{ad}$ with $\bs C \bs p=\bs 0$, we have
 \[\sup_x |\mathbb{P}_{\bs p}\left(Q_N(\bs M)\leq x\right) - \mathbb{P}_{\bs p}\left(Q_N^\epsilon(\bs M)\leq x\right|\bs X)|\stackrel{p}{\to}0.\]
\end{corollary}

The corresponding wild bootstrap version of the ANOVA-type test is given by $\varphi = \mathbbm{1}\{Q_N (\bs M)>c_{Q}^\epsilon(\alpha)\}$, where $c_{Q}^\epsilon(\alpha)$ denotes the conditional $(1-\alpha)$-quantile
of $Q_N^\epsilon(\bs M)$ given the data.

\subsubsection{Multiple Contrast Testing Procedure}\label{sec:MCTPboot}
Similarly to the asymptotic MCTP, a wild bootstrap version of the MCTP can be constructed by means of Theorem \ref{asy1:rm}. For this purpose, a wild bootstrap version of the test statistic is defined by $\bs T_N^\epsilon=(T_1^{p,\epsilon}, \ldots, T_q^{p,\epsilon})'$, where the components
of the vector are given by \[T_\ell^{p,\epsilon}=\frac{\sqrt{N} c_{\ell}' \vwhp^\ep}{\sqrt{\widehat{v}_{\ell\ell}^\epsilon}}, \; \ell=1,\ldots,q.\]
The calculation of $\widehat{v}_{\ell\ell}^\epsilon$ is straightforward using the wild bootstrap version of the empirical covariance matrix $\widehat{\bs V}_N$, which is defined by
$\bs V_N^\epsilon=\bs E_{ad} \cdot \bm{\Sigma}^\epsilon \cdot \bs E_{ad}$.
$\bm{\Sigma}^\epsilon$ is calculated by plugging in the wild bootstrap samples as described in Section \ref{mod}.

\begin{corollary}\label{cor:distr}
If $N/n_i\rightarrow \kappa_i\in(0,\infty)$, the distribution of the wild bootstrap version of the test statistic $\bs T_N^\epsilon$, conditioned on the data, weakly converges to a multivariate normal distribution with mean zero and covariance matrix 
$\bs R$ ($N \to \infty$) in probability.
\end{corollary}

Using Corollary \ref{cor:distr}, the equicoordinate quantile of the normal-${N}(\bs 0,\bs R)$-distribution can be replaced by the corresponding equicoordinate
quantile of the conditional wild bootstrap distribution function of $\bs T_N^\varepsilon$.
Thus, the individual test hypothesis $H_{0,\ell}^p: \bs c'_\ell \bs p=0$ will be rejected if \[|T_\ell^{p}|\geq c^\epsilon(\alpha),\] where $c^\epsilon(\alpha)$ is the conditional $(1-\alpha)$ equicoordinate quantile of $\bs T_N^\epsilon$ given the data. 
Analogously, the global null hypothesis $H_0^p:\bs{Cp}=\bs 0$ will be rejected if \[\max\{|T_1^{p}|, \ldots, |T_q^{p}|\}\geq c^\epsilon(\alpha).\]
Finally, the SCIs for $\bs c'_\ell \bs p$ are given by
\[\left[\bs{c}'_\ell \widehat{\bs{p}}-c^\epsilon(\alpha)\sqrt{\frac{\widehat{v}_{\ell\ell}}{N}};\; \bs{c}'_\ell \widehat{\bs{p}}+c^\epsilon(\alpha)\sqrt{\frac{\widehat{v}_{\ell\ell}}{N}}\right], \; \ell=1,\ldots,q.\]
Additionally to the SCIs for $\bs c'_\ell \bs p$ , a first overview of the contruction of SCIs for ratios is given in the next section.

% \section{Multiple Testing Procedures}\label{stat}% for repeated measures designs} 
% 
% \subsection{The wild bootstrap procedure}\label{sec:MCTPboot}

\subsection{Confidence intervals for ratios}\label{sec:CI}
The construction of SCIs for ratios $\theta_\ell = \bs c_\ell' \bs p / \bs d_\ell' \bs p, \; \ell=1,\dots,q$, where $\bs c_\ell, \bs d_\ell \in \mathbb{R}^{ad}$ are different contrasts ($\bs d_\ell' \bs p\neq 0$),
is based on the mean-based approaches of Dilba et al.~(2004, 2006)\nocite{dilba2004simultaneous}\nocite{dilba2006simultaneous} and Hasler (2009)\nocite{haslerextensions}. In the latter works, SCIs for ratios of the means are constructed, which will
be extended to SCIs for ratios of the relative treatment effect $\bs p$ and the corresponding testing problem
\[H_{0\ell}^\text{ratio}: \theta_\ell = \tau_\ell \quad \text{vs.} \quad H_{1\ell}^\text{ratio}: \theta_\ell < \tau_\ell,  \; \ell=1,\dots,q,\] where $\tau_\ell$ is usually chosen to be 1 for all $\ell=1,\ldots,q$.

Following the ideas of Dilba et al.~(2004, 2006) and Hasler (2009)\nocite{haslerextensions}, the ratio problem $\theta_\ell$ can be expressed by 
the following linear form 
$L_\ell=\left(\theta_\ell \bs d_\ell-\bs c_\ell\right)'\bs p, \; \ell=1,\ldots,q$. %Thus, the confidence intervals can be rested upon a Fieller confidence interval (see Fieller, 1954\nocite{fieller1954some}). 
Then, the vector of test statistics for this ratio problem $\bs T_N^{\text{ratio}}=\left(T_1^{\text{ratio}}, \ldots, T_q^{\text{ratio}}\right)'$ has components 
\[T^{\text{ratio}}_\ell=\frac{\sqrt{N}(\theta_\ell \bs d_\ell-\bs c_\ell)'\left(\widehat{\bs p}-\bs p\right)}{\sqrt{(\theta_\ell \bs d_\ell-\bs c_\ell)'\widehat{\bs V}_N (\theta_\ell \bs d_\ell-\bs c_\ell)}}, \; \ell=1,\ldots,q.\]
Similar to the vector of test statistics $\bs T_N$,  $\bs T_N^{\text{ratio}}$ is asymptotically multivariate normally distributed with mean zero and covariance matrix $\bs S=(s_{\ell m})_{\ell,m=1}^q$.
%with entries
%$s_{\ell m}=\frac{(\theta_\ell \bs d_\ell-\bs c_\ell)' \bs V (\theta_m \bs d_m-\bs c_m)}{\sqrt{(\theta_\ell \bs d_\ell-\bs c_\ell)' \bs V (\theta_\ell \bs d_\ell-\bs c_\ell)} \sqrt{(\theta_m \bs d_m-\bs c_m)' \bs V (\theta_m \bs d_m-\bs c_m)}}$. 
Its distribution can be approximated with the wild bootstrap procedure presented in Section and together with Fieller's Theorem this allows for the construction of wid bootstrap SCIs for the ratios $\theta_\ell, \; \ell=1,\ldots,q$. 
The details and the explicit formula of the SCIs are presented in Appendix \ref{sec:tecdet_ratio}.

\section{Simulation Study} \label{simu}
The behavior of the wild bootstrap procedure for small samples within  an extensive simulations study is examined. Therefore, the maintenance of the nominal type-$I$ error rate of the proposed test procedures is compared. The simulations are conducted with the help of 
\textsc{R} computing environment, version 3.2.3 (R Core Team, 2015\nocite{r2015language}) each with 1,000 simulation runs and 1,000 bootstrap samples. In the first part, the novel multiple testing procedure is compared to the ATS and in the second part
different contrast matrices are compared.
% Our novel wild bootstrap approach (bootMCTP) will be compared to the standard MCTP 
% and the ANOVA-type test statistic (ATS). For the global test, all three different test procedures are compared and for the multiple comparisons, the wild bootstrap procedure is only collated to the standard MCTP.

As in Brunner \& Placzek (2011)\nocite{brunner2011box}, the independent data vectors $\bs X_{ik}, \, i=1,\ldots, a$ and $k=1,\ldots, n_i$ are generated by \[\bs X_{ik}=\sigma_{i} \bs V^{\frac{1}{2}} \bs Z_{ik}+ c_{i}  B_{ik} \bs 1_{d},\]
where $ B_{ik}$ is the effect of the $k$th individual in group $i$ and repeated measure $j$ and $c_{i}$ a scale factor. $\bs Z_{ik}$ generates some error and $\bs V$ denotes the covariance structure of the repeated measures. In this simulation study, 
the vector $\bs Z_{ik}=(Z_{i1k}, \ldots, Z_{idk})'$ is normally distributed with expectation zero and covariance matrix $\bs I_{d}$ and $\bs B_{i}=(B_{i1}, \ldots, B_{in_i})'$ is chosen to be normally distributed with expectation zero and 
covariance matrix $\bs I_{n_i}$.
Furthermore, three different covariance structures are taken into account:
\begin{itemize}
 \item $\bs V=\bs I_{d}$ (compound symmetry structure, CS),
 \item $\bs V=(v_{\ell m})_{\ell,m=1}^d=\rho^{|\ell-m|}$, where $\rho \in (0,1)$ (autoregressive structure, AR($\rho$)),
 \item $\bs V=(v_{\ell m})_{\ell,m=1}^d=d-|\ell-m|$ (Toeplitz structure, TPL).
\end{itemize}
Regarding this simulations, the constant $\rho$ of the autoregressive covariance structure is determined to be 0.6.

Two different balanced, homoscedastic designs are simulated. First, a design including three different levels of a treatment over three different time points for each individual is examined. Hereafter, this design is called Setting 1. Similar to Setting 1,
a model regarding two different treatment groups measured at four different time points are conducted, which is denoted by Setting 2 in the following. In the following, six different sample sizes (5, 10, 15, 20, 25, 30) are compared and all 
three covariance structures for the repeated measures introduced above are considered. 

\subsection{Comparisons with global testing methods}

%The simulation study is divided into two parts. 
In the first part, the ANOVA-type test statistic (ATS) is compared to the standard and the wild bootstrap MCTP. Therefore, it only makes sense to use an adequate centring
matrix $\bs P$ (with the right dimensions) as presented at the end of Section 1 as a contrast matrix. 
% In the second part, the two multiple contrast tests should be examined. Thus, the ATS is not taken into consideration, since the ATS is not an adequate testing procedure for arbitrary
% contrast matrix as e.g. the Tukey-, Dunnett- or Average-type contrast matrix.

\begin{figure}[!ht]
\centering
 \includegraphics[width=\textwidth]{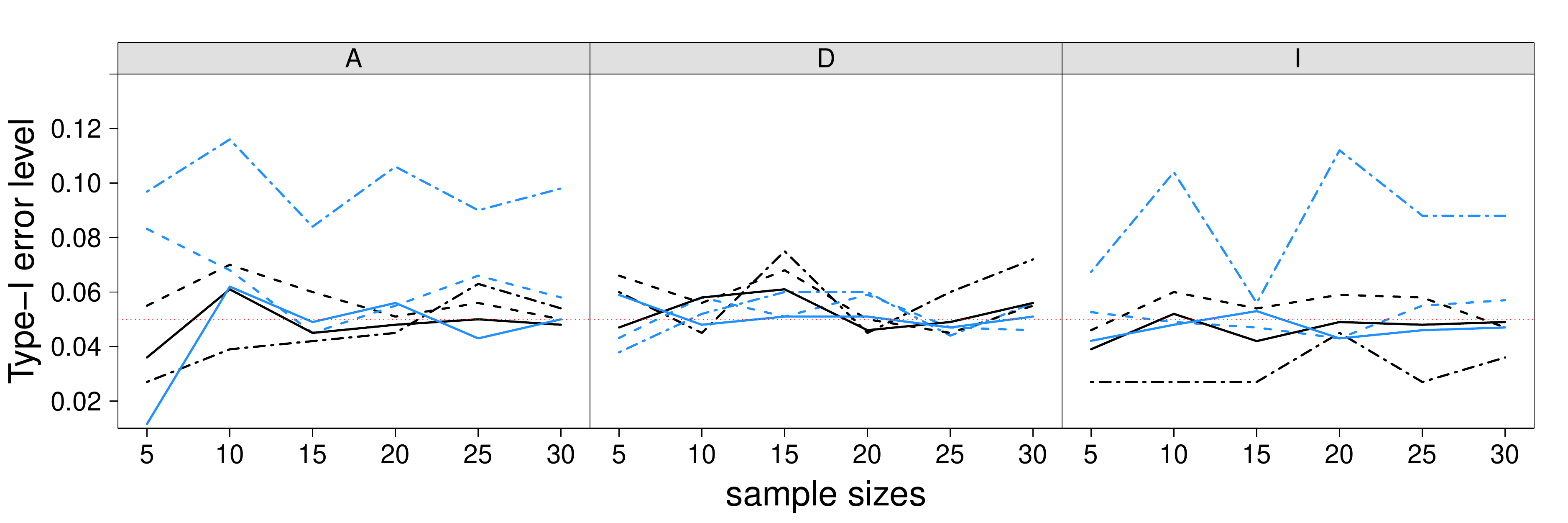}
 \caption{\label{fig:CS}Type-$I$ error rates for Setting 1 (black) and Setting 2 (blue) with a compound symmetry covariance structure  for three different tests, namely MCTP (dotted-dashed), bootMCTP (solid) and ATS (dashed).}
\end{figure}

\begin{figure}[!ht]
\centering
 \includegraphics[width=\textwidth]{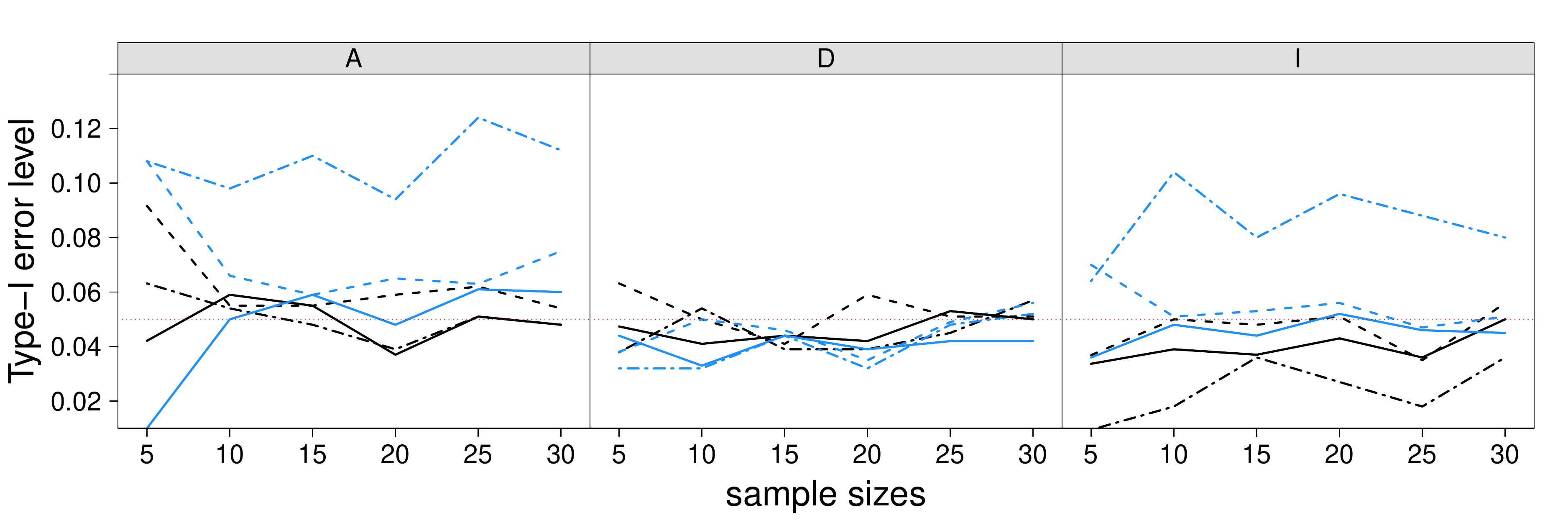}
 \caption{\label{fig:AR}Type-$I$ error rates for Setting 1 (black) and Setting 2 (blue) with an autoregressive covariance structure  for three different tests, namely MCTP (dotted-dashed), bootMCTP (solid) and ATS (dashed).}
\end{figure}

\begin{figure}[!ht]
\centering
 \includegraphics[width=\textwidth]{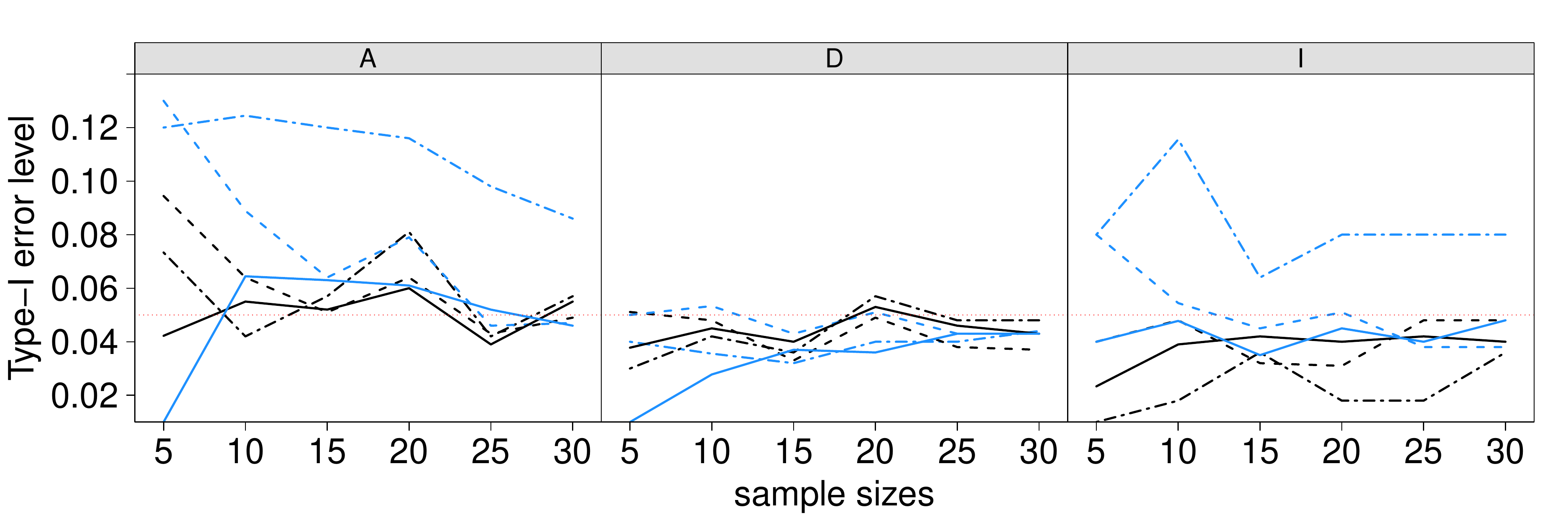}
 \caption{\label{fig:TPL}Type-$I$ error rates for Setting 1 (black) and Setting 2 (blue) with a Toeplitz covariance structure for three different tests, namely MCTP (dotted-dashed), bootMCTP (solid) and ATS (dashed).}
\end{figure}

The results for the compound symmetry structure are summarized in Figure~\ref{fig:CS}, Figure~\ref{fig:AR} visualizes the results for the autoregressive structure and Figure~\ref{fig:TPL} for the 
Toeplitz structure.
The difference between the three different covariance structures is negligible.  Even to find the best testing procedure is quite difficult since all three approaches have their advantages and disadvantages. In some cases, the ATS yield better results 
than the bootstrap version of the MCTP. But note, that the ATS is not an adequate testing procedure for multiple comparisons and therefore, only applicable for global testing problems. Furthermore, the wild bootstrap multiple testing procedure
shows pretty good results in controlling the type-$I$ errors for rising sample sizes, when compared to the ATS and to the standard MCTP. Especially when comparing the standard and the wild bootstrap based MCTP, there are 
some cases (``no main effect $A$'' and ``no interaction effect'') the standard MCTP shows a very liberal behavior, whereas the wild bootstrap MCTP exhibits accurate type-$I$ error level control in almost all cases.

\subsection{Investigating the impact of different contrast matrices}
In this subsection, the novel multiple testing procedures are compared with regard to different contrast matrices. Thus, the ATS is not taken into consideration, because the tests are not comparable in this situation. 
\begin{figure}[!ht]
\centering
 \includegraphics[width=\textwidth]{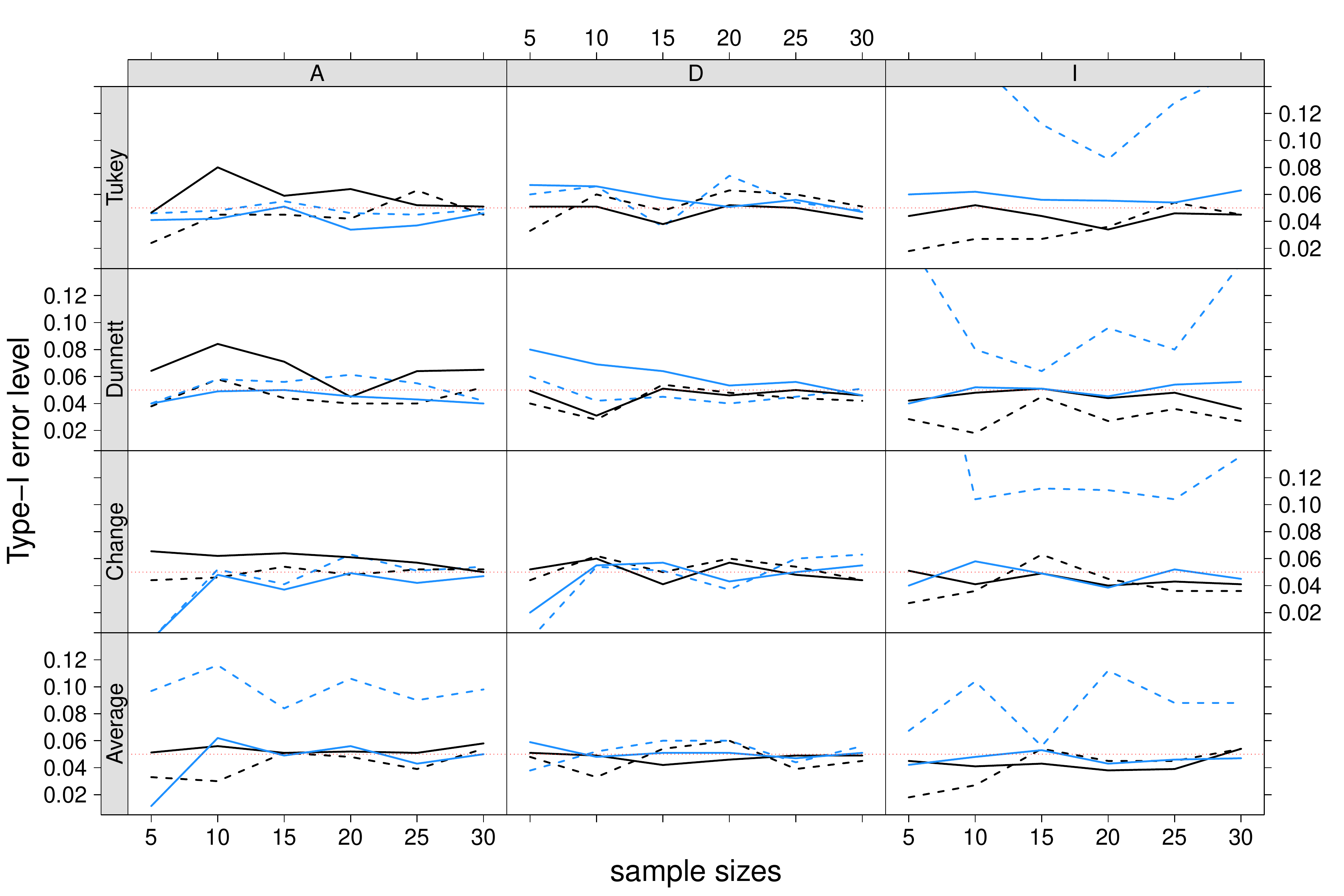}
 \caption{\label{fig:erg_CS_mult}Type-$I$ error rates for Setting 1 (black) and Setting 2 (blue) with a compound symmetry covariance structure for four contrast matrices and two tests, namely the MCTP (dashed) and the bootMCTP (solid).}
\end{figure}
\begin{figure}[!ht]
\centering
 \includegraphics[width=\textwidth]{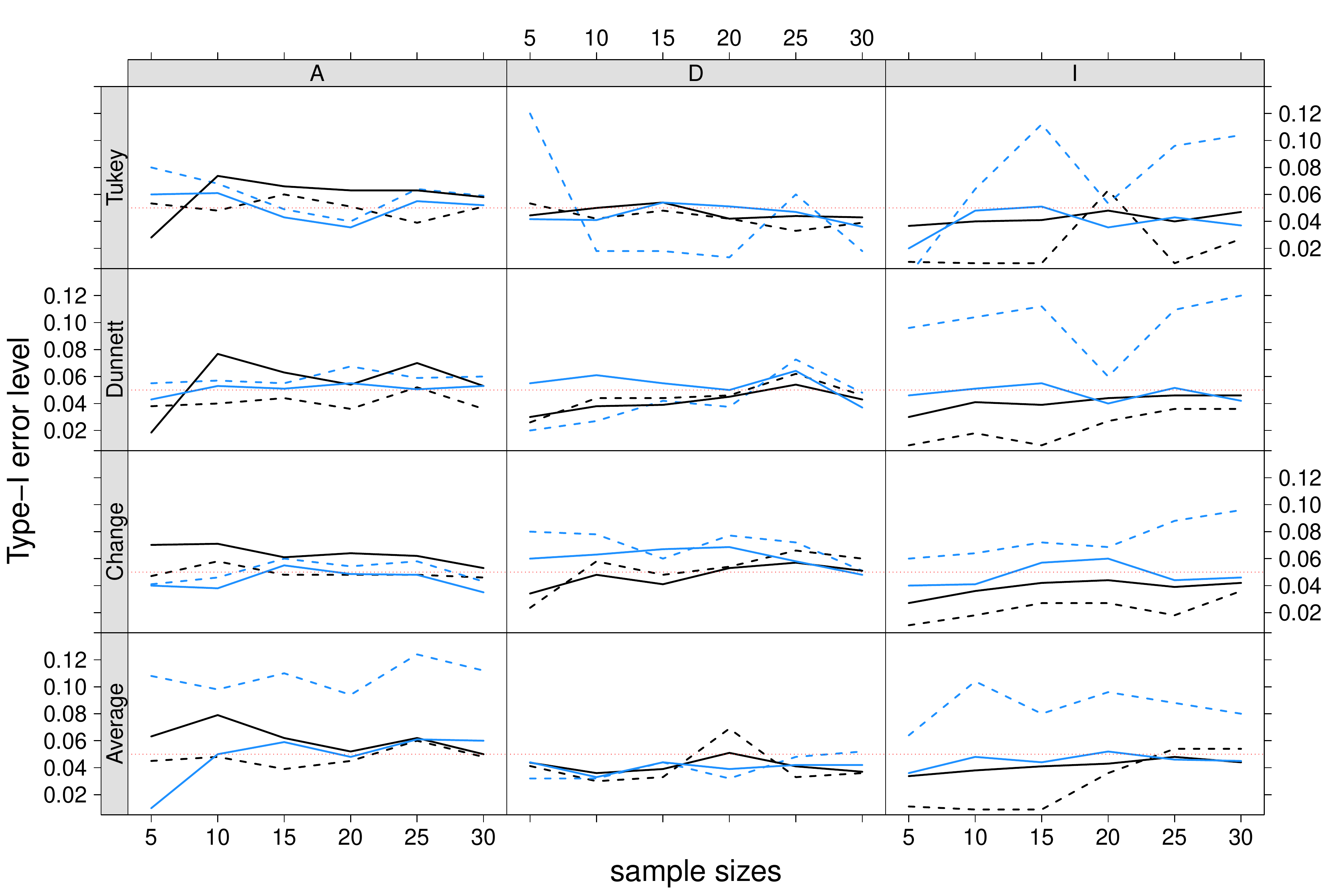}
 \caption{\label{fig:erg_AR_mult}Type-$I$ error rates for Setting 1 (black) and Setting 2 (blue) with an autoregressive covariance structure for four contrast matrices and two tests, namely the MCTP (dashed) and the bootMCTP (solid).}
\end{figure}
\begin{figure}[!ht]
\centering
\includegraphics[width=\textwidth]{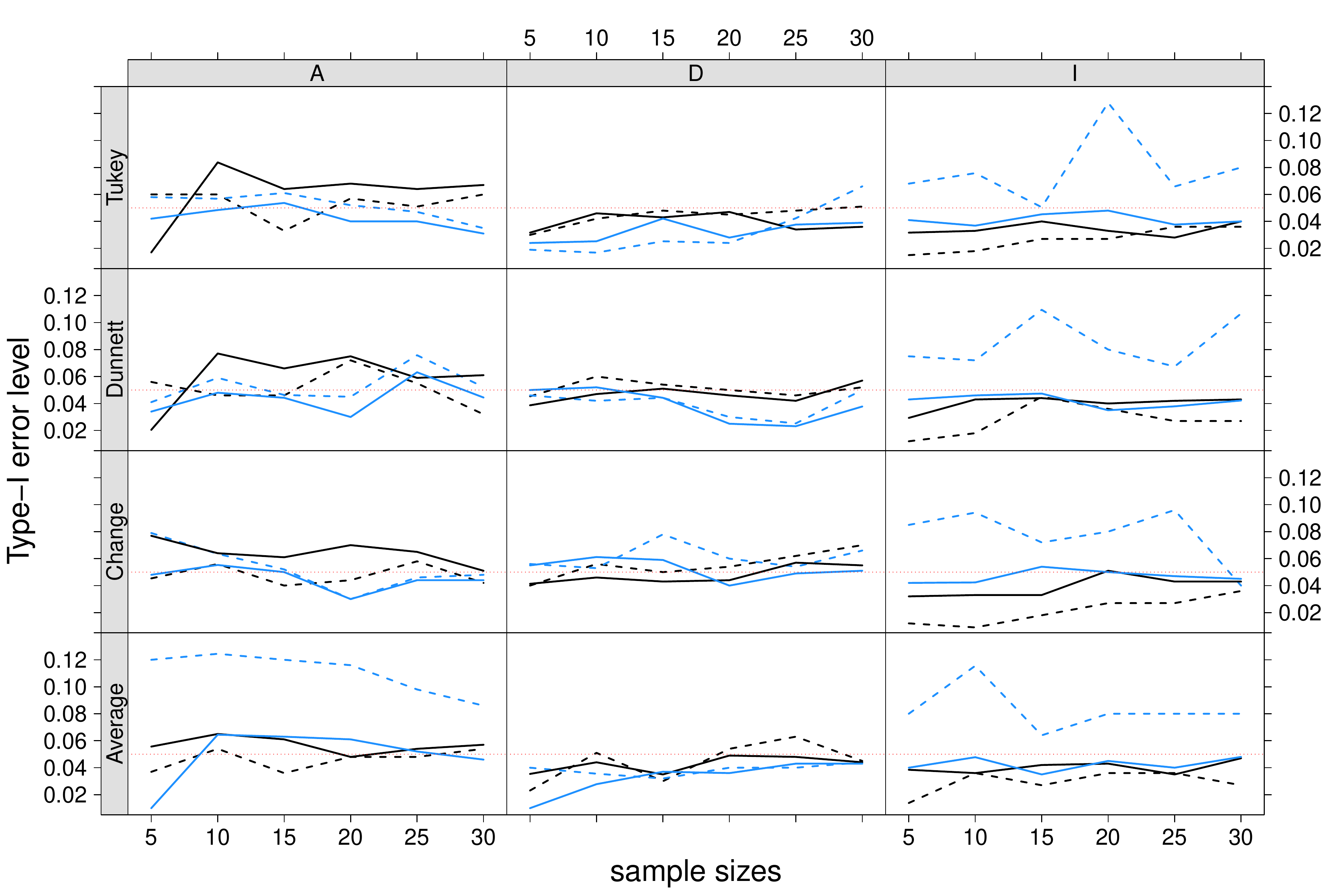}
 \caption{\label{fig:erg_TPL_mult}Type-$I$ error rates for Setting 1 (black) and Setting 2 (blue) with a Toeplitz covariance structure for four contrast matrices and two tests, namely the MCTP (dashed) and the bootMCTP (solid).}
\end{figure}
%In the second part of this simulation study, different contrast matrices are compared. 
Here, four different contrast matrix are taken into account; namely for all-pairs (Tukey), average, many-to-one (Dunnett) and changepoint comparisons. All corresponding 
contrast matrices are summarized in Appendix~\ref{sec:contmat}. Again, two different settings, three different covariance structures, and six different sample sizes are examined. The results are summarized in 
Figures~\ref{fig:erg_CS_mult}-\ref{fig:erg_TPL_mult} regarding the compound symmetry structure in the first, the autoregressive structure in the second and the Toeplitz structure in the third of this three figures.

Even in case of a comparison between the standard and the wild bootstrap based MCTP, the results between the different covariance structures and the different contrast matrices are quite similar. Some distinctions can be made for the Average-type contrast matrix
and the interaction. Regarding the Average-type contrast matrix and the results for \textit{no main effect $A$}, the standard MCTP shows very liberal behavior in Setting 2. In case of the \textit{no interaction effect}, the standard MCTP shows a liberal
behavior for Setting 2 and tends to conservative results in Setting 1, whereas the wild bootstrap MCTP controls the type-$I$ error very accurately.

%The behavior of the standard MCTP and the wild bootstrap version of the MCTP in the case of multiple contrast test is examined more precisely in the next section, where the MCTP is applied to two different data examples.

\section{Application to empirical data}\label{exam}

Now, the theoretical statements made above are applied to empirical data. First, a dataset included in the \textsc{R}-package \texttt{nparLD} is examined and afterwards data from the Institute of Clinical and Biological Psychology at Ulm University 
is analyzed.

\subsection{Shoulder tip pain study}
The following dataset (\texttt{shoulder}) is obtained from the \textsc{R}-package \texttt{nparLD}. The dataset was also studied by Lumley (1996)\nocite{lumley1996generalized}.

\begin{figure}[ht]
\centering
\includegraphics[width=0.9\textwidth]{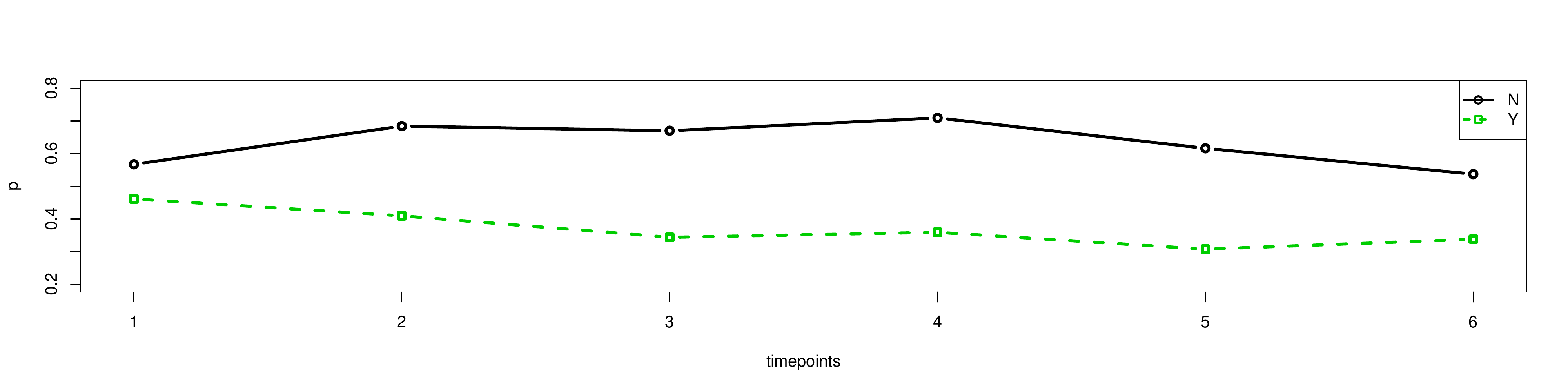}
\caption{\label{fig:exam}The relative effects of two different treatment (``Y'' and ``N'') and six different time points of the \texttt{shoulder} dataset.}
\end{figure}

In this study 
the shoulder pain level of 41 patients after a laparoscopic surgery in the abdomen was examined. During such a surgery, the surgeon fills the abdominal part of the patient with air to have a better view of the body. After this laparoscopic surgery, 
the air was removed out of the abdomen by using a specific suction procedure. A random subsample of 22 patients (``Y'') was treated with this special suction method. In the other subsample of 19 patients (``N'') the air was left in the abdomen. 
The patients were asked for their pain score two times a day (morning and evening) for the first 
three days after the surgery and this score has five different levels from 1 (= low) to 5 (= high). The relative effects of the data for the six different time points are given in Figure~\ref{fig:exam}.

In the following, we like to work out the difference between the different groups. Therefore, we apply a Tukey-type contrast matrix to make all-pairs comparisons. Using this results, we like to find out which of the timepoint differ from each 
other. The results are summarized in Table~\ref{tab:tukey_pain}. 

\begin{table}[ht]
 \centering
 \footnotesize
 \caption{\label{tab:tukey_pain}Many-to-one comparison of the shoulder tip pain study for the standard MCTP in the middle part and for the wild bootstrap MCTP in the right part of the table, significant values are printed in bold.}
 \begin{tabular}{lr|crc|cc}
  \toprule
  Comparison & $\widehat{p}_{\cdot i}-\widehat{p}_{\cdot i'}$ & 95\%-CI& $t$-value &$p$-value & 95\%-CI (wb) & $p$-value (wb) \\
  \midrule
  timepoint 2 vs. timepoint 1 & 0.033 & $\left[-0.091 ; \hspace{0.8em} 0.155\right]$ & 0.781 & 0.963 & $\left[-0.045 ; \hspace{0.8em}0.111\right]$ & 0.443\\
  timepoint 3 vs. timepoint 1 & -0.007 & $\left[-0.124 ; \hspace{0.8em}0.110\right]$ & -0.185& 0.999 & $\left[-0.081 ; \hspace{0.8em}0.067\right]$ & 0.858\\
  timepoint 4 vs. timepoint 1 & 0.020 & $\left[-0.095 ; \hspace{0.8em}0.134\right]$ & 0.515 & 0.994 & $\left[-0.058 ; \hspace{0.8em}0.098\right]$ & 0.624\\
  timepoint 5 vs. timepoint 1 & -0.052 & $\left[-0.171 ; \hspace{0.8em}0.068\right]$ & -1.282 & 0.768 & $\left[-0.127 ; \hspace{0.8em}0.023\right]$ & 0.193\\
  timepoint 6 vs. timepoint 1 & -0.077 & $\left[-0.186 ; \hspace{0.8em}0.035\right]$ &-2.028 &0.319 & $\left[-0.145; -0.008\right]$&\textbf{0.035}\\
  timepoint 3 vs. timepoint 2 & -0.040 & $\left[-0.111 ; \hspace{0.8em}0.031\right]$ &-1.657 &0.536 & $\left[-0.087 ; \hspace{0.8em} 0.007\right]$ &0.091\\
  timepoint 4 vs. timepoint 2 & -0.013 & $\left[-0.092 ; \hspace{0.8em}0.067\right]$ &-0.469 & 0.996 & $\left[-0.064;  \hspace{0.8em}0.038\right]$&0.629\\
  timepoint 5 vs. timepoint 2 & -0.085 & $\left[-0.184 ; \hspace{0.8em} 0.016\right]$ &-2.480 & 0.141 & $\left[-0.151; -0.018\right]$ &\textbf{0.007}\\
  timepoint 6 vs. timepoint 2 & -0.109 & $\left[-0.213 ; -0.003\right]$ & -3.005 & \textbf{0.045} & $\left[-0.178; -0.038\right]$ & \textbf{0.003}\\
  timepoint 4 vs. timepoint 3 & 0.027 & $\left[-0.043 ; \hspace{0.8em}0.097\right]$ & 1.150 & 0.838 & $\left[-0.021; \hspace{0.8em}  0.075\right]$ &0.253\\
  timepoint 5 vs. timepoint 3 & -0.045 & $\left[-0.139 ; \hspace{0.8em}0.050\right]$ & -1.404 & 0.696 & $\left[ -0.105;  \hspace{0.8em}0.015\right]$ & 0.148\\
  timepoint 6 vs. timepoint 3 & -0.069 & $\left[-0.166 ; \hspace{0.8em}0.029\right]$ & -2.075 & 0.296 & $\left[-0.133 ;-0.004\right]$ & \textbf{0.031}\\
  timepoint 5 vs. timepoint 4 & -0.072 & $\left[-0.149 ; \hspace{0.8em}0.006\right]$ & -2.731 & 0.084 & $\left[-0.120; -0.024\right]$& \textbf{0.005}\\
  timepoint 6 vs. timepoint 4 & -0.097 & $\left[-0.176 ; -0.015\right]$ & -3.479 & \textbf{0.014} & $\left[-0.152; -0.041\right]$& \textbf{0.000}\\
  timepoint 6 vs. timepoint 5 & -0.024 & $\left[-0.072 ; \hspace{0.8em}0.024\right]$ & -1.487 & 0.644 & $\left[-0.056;  \hspace{0.8em}0.008\right]$& 0.145\\
   \midrule
  Y vs. N & -0.261 & $\left[-0.365 ; -0.138\right]$ & -4.525 & \textbf{$<$0.001} & $\left[-0.365; -0.139\right]$ &\textbf{$<$0.001}\\
  \bottomrule
 \end{tabular}
\end{table}

% \begin{table}[!ht]
% \centering
%  % \footnotesize
%   \begin{tabular}{ll|ccc}
%   \toprule
%      & & treatment & time &interaction\\
%       \midrule
%   \multirow{3}{*}{all-pairs} & ATS & \textbf{$<$0.001} & \textbf{0.031} & \textbf{0.011}\\
%  & MCTP & \textbf{$<$0.001} & \textbf{0.014} & 0.055 \\
%  & wildMCTP & \textbf{$<$0.001} & \textbf{$<$0.001}  &\textbf{0.013} \\
%   \midrule
%     \multirow{3}{*}{many-to-one} & ATS &  \textbf{$<$0.001} & \textbf{0.031} & \textbf{0.011}\\
%  & MCTP & \textbf{$<$0.001} & 0.131 &  0.055\\
%  & wildMCTP & \textbf{$<$0.001} & \textbf{0.037} & \textbf{0.013}\\
% %   \midrule
% %     \multirow{3}{*}{changepoint} & ATS &\textbf{$<$0.001} &  \textbf{0.031} & \textbf{0.011}\\
% %  & MCTP & \textbf{$<$0.001} & \textbf{0.007} & 0.055 \\
% %  & wildMCTP & \textbf{$<$0.001}  & \textbf{0.001} &  \textbf{0.013}\\	
%   \bottomrule
%   \end{tabular}
%   \caption{\label{tab:exam}p-values for the global test of the shoulder tip pain study with 999 bootstrap iterations, significant values are printed in bold.}
% \end{table}

First of all, one can see that in almost all cases the CIs of the wild bootstrap MCTP show much shorter widths than the CIs of the standard multiple testing methods. Thus, the wild bootstrap version of the MCTP yields more 
statistically significant values than the standard version of the MCTP. All of the significant values in the standard MCTP are also significant for the wild bootstrap MCTP. Moreover, the test procedures are consonant and coherent since all $p$-values
correspond to the CIs, that is, CIs not including zero lead to a significant $p$-value. With references to the shoulder tip pain study, the standard MCTP yields two pairs of groups,
which are statistically different, namely timepoint 2 to timepoint 6 ($p=0.045$) and timepoint 4 to timepoint 6 ($p=0.014$). Beyond these two pairs, the wild bootstrap approach detects four other pairs which vary from one to the other group. These are
timepoint 1 to timepoint 6 ($p=0.035$), timepoint 2 to timepoint 5 ($p=0.007$), timepoint 3 to timepoint 6 ($p=0.031$) and timepoint 4 to timepoint 5 ($p=0.005$).

\begin{figure}[ht]
\centering
\includegraphics[height=0.3\textheight]{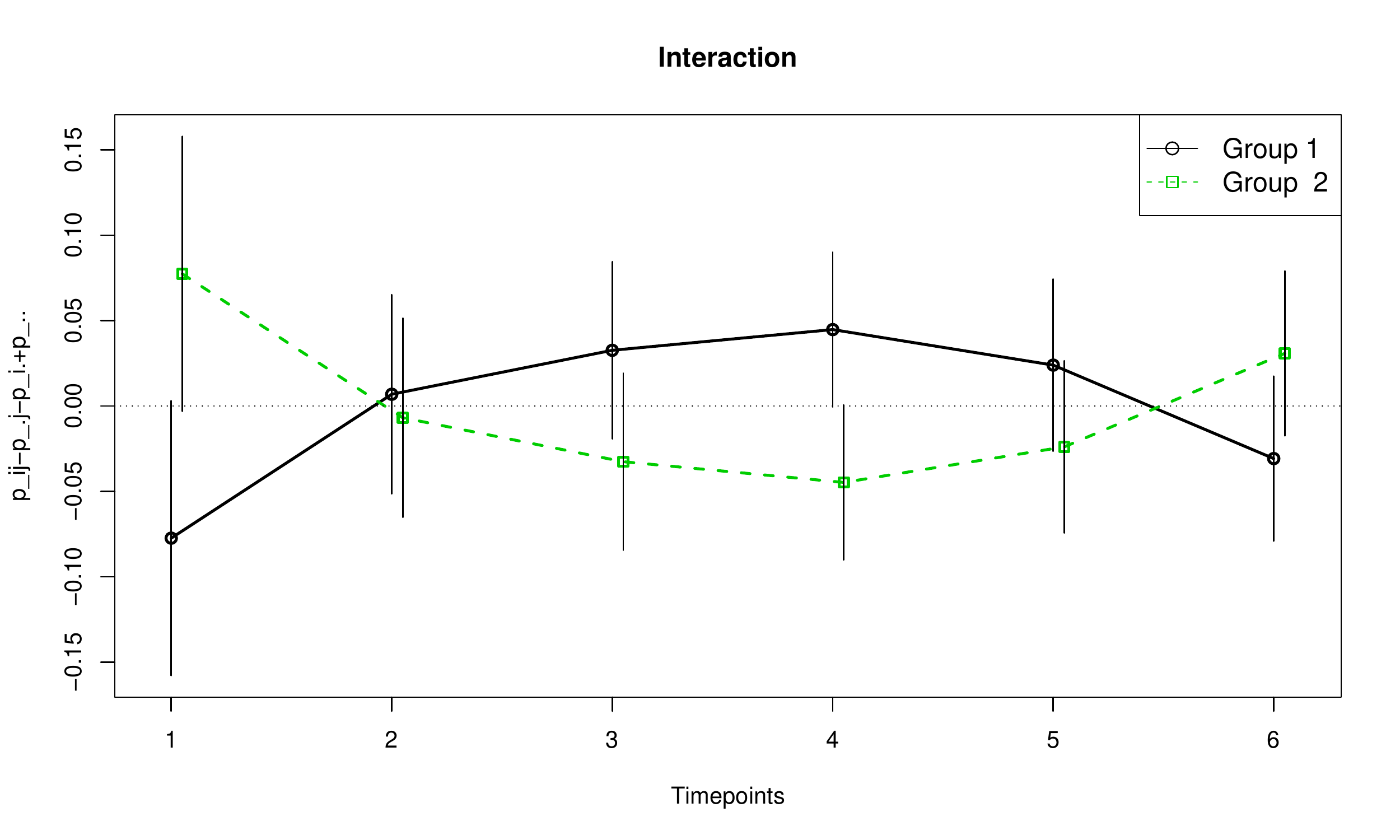}
\caption{\label{fig:pain_int}Interaction plot for the shoulder tip pain study.}
\end{figure}

Moreover, an interaction plot for the shoulder tip pain study is given in Figure~\ref{fig:pain_int}. Generally, in case of a significant interaction effect, one may be interested in the question: ``Which profile of which treatment group differs?''. 
And the courses of such an interaction plot as in Figure~\ref{fig:pain_int} may be the first hint. In a second step, pairwise comparisons should be conducted, since only looking at a plot did not verify statistical significance.
In case of the shoulder tip pain study, the wild bootstrap version of the MCTP and the ATS yield a significant result (wildMCTP: $p=0.013$, ATS: $p=0.011$), whereas the $p$-value of the standard MCTP is 0.055. As one can the in Figure~\ref{fig:pain_int},
the courses of the two treatment groups are nearly mirrored at the x-axis and thus, are not equal. An all-group comparison, in this case only a comparison between the two treatment groups, yields the following result: The calculated $p$-values for the
standard and the wild bootstrap MCTP are $<0.001$ and thus, this inference procedure confirms the latter assumption.

\subsection{Childhood maltreatment study}

\begin{figure}[ht]
\centering
\includegraphics[width=0.9\textwidth]{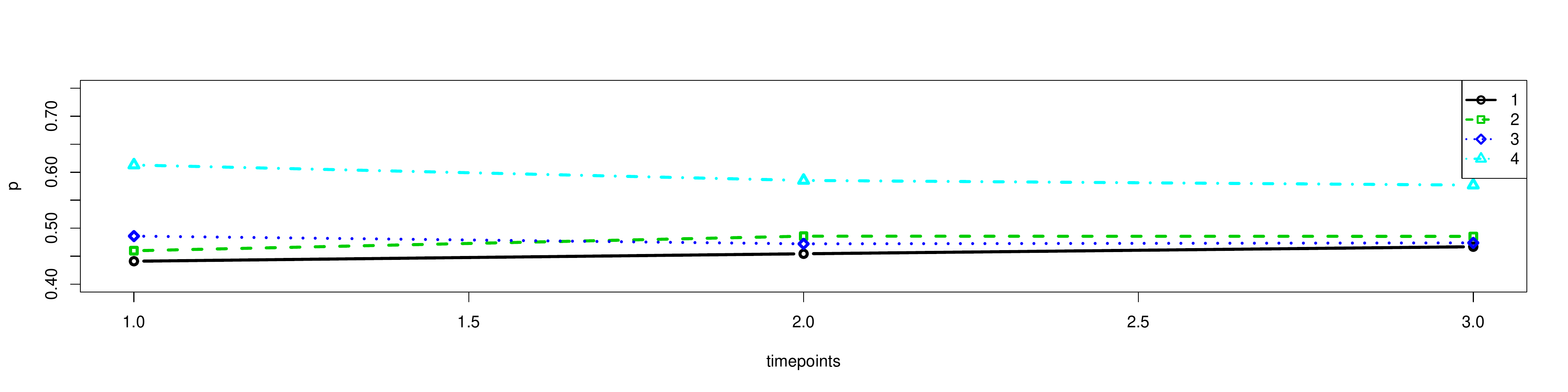}
\caption{\label{fig:exam2}The relative effects of women with different CTQ values (categories 1-4) and three different time points.}
\end{figure}

In this section, a partial dataset of a recent study from the Institute of Clinical and Biological Psychology at Ulm University on childhood maltreatment, postnatal distress and the role of social support is examined. 
Using the dataset, we like to 
determine the influence of the childhood maltreatment on the postnatal distress of the mothers. The postnatal distress was measured three times, three ($t_1$), six ($t_2$) and nine month ($t_3$) postpartum. The postnatal psychological distress 
(dependent variable) was estimated by a combined score of the Perceived Stress Scale (PSS4) and the Hospital Anxiety and Depression Scale (HADS). For the postnatal psychological distress, sum scores of the PSS4 and HADS scales were standardized and 
added together. The maltreatment experience  such as emotional, physical and sexual abuse as well as emotional and physical neglect in the mother’s own childhood was measured with the Childhood Trauma Questionnaire (CTQ). Using quantiles, the CTQ sum 
score was categorized, resulting in the following value ranges: 25-29 (first quartile), 30-32 (median), 33-37 (third quartile) and 38-103 (maximal value). 
The relative effects of the data for the three different time points are given in Figure~\ref{fig:exam2}. 

Figure~\ref{fig:exam2} shows that the estimates of the relative effects of the first three groups are nearly the same. Only the values of the relative effects of the group including the highest CTQ
values (group 4) are higher compared to the other groups.

We like to work out the difference between the different CTQ groups and the different measurements. Therefore all-pairs comparisons in both factors are conducted. The results regarding the CTQ categories and the three measurement time points are presented
in Table~\ref{tab:tukey_CM}.

\begin{table}[ht]
 \centering
 \footnotesize
 \caption{\label{tab:tukey_CM}Many-to-one comparison of the childhood maltreatment study for the standard MCTP in the middle part and for the wild bootstrap MCTP in the right part of the table.}
 \begin{tabular}{lr|crc|cc}
  \toprule
  Comparison & $\widehat{p}_{\cdot i}-\widehat{p}_{\cdot i'}$ & 95\%-CI& $t$-value &$p$-value & 95\%-CI (wb) & $p$-value (wb) \\
   \midrule
  category 2 vs. category 1 & 0.023 & $\left[-0.261 ; 0.303\right]$ & 0.229 & 0.995 & $\left[-0.174 ; 0.218\right]$ & 0.821\\
  category 3 vs. category 1 & 0.023 & $\left[-0.237 ; 0.280\right]$ & 0.253 & 0.994 & $\left[-0.164 ; 0.208\right]$ & 0.804\\
  category 4 vs. category 1 & 0.138 & $\left[-0.213 ; 0.455\right]$ & 1.112 & 0.685 & $\left[-0.109 ; 0.367\right]$ & 0.271\\
  category 3 vs. category 2 & 0.000 & $\left[-0.290 ; 0.291\right]$ & 0.002 & 1.000 & $\left[-0.196 ; 0.196\right]$ & 0.997\\
  category 4 vs. category 2 & 0.115 & $\left[-0.253 ; 0.453\right]$ & 0.881 & 0.812 & $\left[-0.126 ; 0.342\right]$ & 0.400\\
  category 4 vs. category 3 & 0.114 & $\left[-0.241 ; 0.442\right]$ & 0.913 & 0.796 & $\left[-0.117 ; 0.332\right]$ & 0.360\\
  \midrule
  timepoint 2 vs. timepoint 1 & -0.001 & $\left[-0.109 ; 0.108\right]$ & -0.011 & 1.000 & $\left[-0.088 ; 0.086\right]$ & 0.988\\
  timepoint 3 vs. timepoint 1 & 0.001 & $\left[-0.108 ; 0.110\right]$ & 0.021 & 1.000 & $\left[-0.086 ; 0.088\right]$ & 0.987\\
  timepoint 3 vs. timepoint 2 & 0.001 & $\left[-0.021 ; 0.023\right]$ & 0.159 & 0.985 & $\left[-0.015 ; 0.017\right]$ & 0.885\\
  \bottomrule
 \end{tabular}
\end{table}

Regarding Figure~\ref{fig:exam2}, a difference between the time point is not expected since the four presented relative effect courses are nearly straight lines. Nevertheless, a difference of CTQ category 4 (high CTQ values) to all other groups could be expected. However,
no significant values -- neither in the whole-plot (CTQ category) nor in the sub-plot factor (time) -- were calculated. Since the estimated values of the differences of the treatment effects in the time factor are near to zero, all calculated
CIs include this value and no significant $p$-value can be computed.  But again, in all cases, the CIs regarding the wild bootstrap approach are much shorter than in case of the 
standard procedure.

\begin{figure}[ht]
\centering
\includegraphics[height=0.3\textheight]{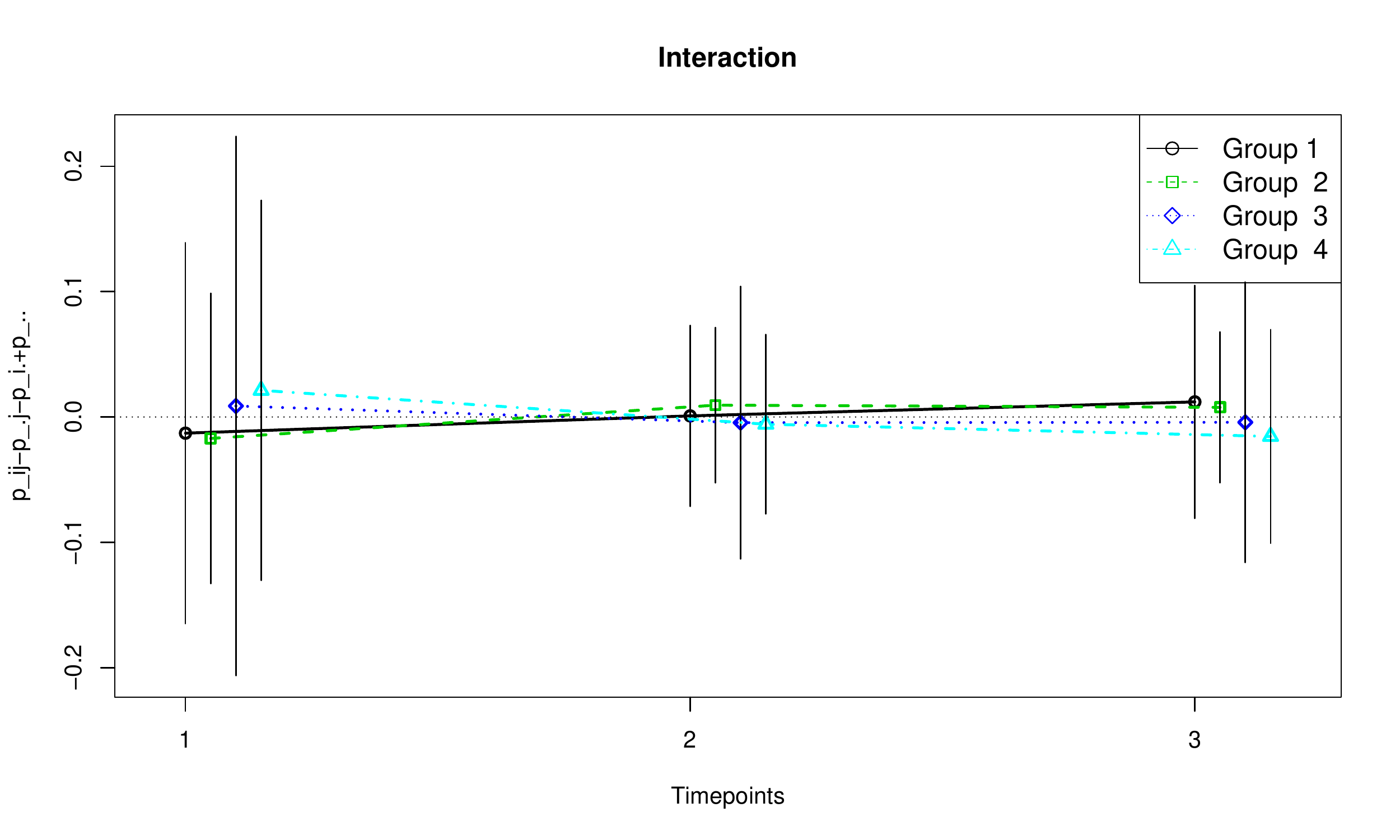}
\caption{\label{fig:cm_int}Interaction plot for the childhood maltreatment study.}
\end{figure}

Additionally, an interaction plot is given in Figure~\ref{fig:cm_int}. In case of the childhood maltreatment study, the $p$-value for the global null hypothesis of no interaction effect was not significant. The standard MCTP yields a 
$p$-value of 0.9771, the $p$-value of the ATS was 0.1290, and for the wild bootstrap version of the MCTP, a $p$-value of 0.5856 was computed. All values indicate no different course of the interaction profiles and thus, Figure~\ref{fig:cm_int} does, since the profiles are nearly the same.

\section{Discussion and Conclusion} \label{con}
Konietschke et al.~(2012)\nocite{konietschke2012rank} already developed multiple, rank-based contrast tests for one-way layouts and also derive SCIs. Additionally, the recent work of Brunner et al.~(2018)\nocite{brunner2016rank}
introduces results for adequate effect measures in repeated measures. Extending both approaches leads to a novel asymptotically exact multiple testing procedure for general factorial split-plot designs. Additionally, the technical details of a wild
bootstrap approach for the test statistic were given. Moreover, SCIs for split-plot designs were introduced and also novel SCIs for ratios of different contrast were constructed.

Since the dependence structure within the data becomes much more complex when repeated measures designs are present, classical inference methods to evaluate such data show their limits. Furthermore, the used effect sizes may not be appropriate. Thus, a 
multiple rank-based testing procedure for split-plot designs and also a wild bootstrap version to obtain asymptotically exact results are studied in this work. The asymptotical exactness of the global testing procedure is an improvement over the
classical ATS. Simulation studies show that the novel MCTP controls the type-$I$ error level quite accurately. And the applicability of the proposed methods was indicated by the analysis of two different real data examples.

In future works, we like to extend the proposed methods to clustered data, which are often present in clinical studies or other practical applications.

\newpage
\section*{Acknowledgements}

The work of Maria Umlauft and Markus Pauly was supported by the German Research Foundation project DFG-PA 2409/3-1.

\nocite{dobler2017nonparametric}

\bibliography{biblio}

\newpage
\appendix

\section{Contrast Matrices}\label{sec:contmat}
The Tukey-type contrast matrix is given by
\[\bs C=\begin{pmatrix}
     -1 & 1 & 0 & 0 & \ldots & 0 & 0 \\
     -1 & 0 & 1 & 0 & \ldots & 0 & 0 \\
     \vdots& \vdots & \vdots & \vdots & \vdots & \vdots & \vdots\\
     -1 & 0 & 0 & 0 & \ldots & 0 & 1 \\
     0 & -1 & 1 & 0 & \ldots & 0 & 0 \\
     0 & -1 & 0 & 1 & \ldots & 0 & 0 \\
     \vdots& \vdots & \vdots & \vdots & \vdots & \vdots & \vdots\\
     0 & 0 & 0 & 0 & \ldots & -1 & 1 \\
    \end{pmatrix},
\]
the Dunnett-type contrast matrix by
\[\bs C=\begin{pmatrix}
         -1 & 1 & 0 & 0 & \ldots & 0 & 0\\
         -1 & 0 & 1 & 0 & \ldots & 0 & 0\\
         \vdots& \vdots & \vdots & \vdots & \vdots & \vdots & \vdots\\
         -1 & 0 & 0 & 0 & \ldots & 0 & 1\\
        \end{pmatrix},
\]
the Average-type contrast matrix by
\[\bs C = \begin{pmatrix}
           1 & -\frac{1}{a-1} & -\frac{1}{a-1} & \ldots & -\frac{1}{a-1} \\
           -\frac{1}{a-1} & 1 & -\frac{1}{a-1} & \ldots & -\frac{1}{a-1} \\
	   \vdots & \vdots & \ddots & \vdots &\vdots \\
	   \vdots & \vdots & \vdots & \ddots &\vdots \\
           -\frac{1}{a-1} & -\frac{1}{a-1} & -\frac{1}{a-1} & \ldots & 1 \\
          \end{pmatrix},
\]
and finally, the matrix for the changepoint comparisons
\[\bs C= \begin{pmatrix}
          -1 & \frac{n_2}{\sum_{i=2}^a n_i} & \frac{n_3}{\sum_{i=2}^a n_i} & \ldots & \frac{n_{a-1}}{\sum_{i=2}^a n_i} & \frac{n_a}{\sum_{i=2}^a n_i} \\ 
          -\frac{n_1}{\sum_{i=1}^2 n_i} & -\frac{n_2}{\sum_{i=1}^2 n_i} & \frac{n_3}{\sum_{i=2}^a n_i} & \ldots & \frac{n_{a-1}}{\sum_{i=2}^a n_i} & \frac{n_a}{\sum_{i=2}^a n_i} \\
          \vdots & \vdots & \ddots & \vdots & \vdots & \vdots \\
          \vdots & \vdots & \vdots & \ddots & \vdots & \vdots \\
          -\frac{n_1}{\sum_{i=1}^{a-1} n_i} & -\frac{n_2}{\sum_{i=1}^{a-1} n_i} & -\frac{n_3}{\sum_{i=1}^{a-1} n_i} & \ldots & -\frac{n_{a-1}}{\sum_{i=1}^{a-1} n_i} & 1 \\
         \end{pmatrix}.
\]
The matrices presented in this section are only a small selection of possible contrast matrices. Because of comparability to extisting simulation studies (e.g. see Konietschke et al., 2012\nocite{konietschke2012rank}) these contrast
matrices were chosen for our simulations.

\section{Technical Details}\label{sec:tecdet}
\subsection{Explicit formulas of the variances and covariances}\label{sec:tecdet_cov}
Here, the sophisticated results of the variance $\sigma_{rs}(p,q,p',q')$ and the covariance $\sigma_{rs,il}(p,q,p',q')$ defined in Section \ref{mod} are given. 

Assuming $X_{ijk}$ and $X_{i'j'k'}$ are independent for $i \neq i'$ or $k \neq k'$, the entries of the asymptotic covariance matrix are given by
\begin{align}\label{equ:asycovrs}
\begin{split}
 & \frac{\sigma_{rs}(p,q,p',q')}{N}\\& \qquad=\left\{\begin{tabular}{p{5.3cm}p{5cm}} $\tau_r^{(s,s)}(p,q,p',q')$, & $r \notin \{p,p'\}, \; p \neq p'$,\\[1em]
						   $\tau_r^{(s,s)}(p,q,p,q') + \tau_p^{(q,q')}(r,s,r,s)$, & $r \notin \{p,p'\}, \; p = p'$,\\[1em]
						   $\tau_r^{(s,s)}(r,q,p',q') - \tau_r^{(q,s)}(r,s,p',q')$,& $r = p, \; p' \neq p', \; q \neq s$,\\[1em]
						   $\tau_r^{(s,s)}(p,q,r,q') - \tau_r^{(s,q')}(p,q,r,s)$,& $r=p, \; p' \neq p, \; q' \neq s $,\\[1em]
						   $\tau_r^{(s,s)}(p,q,r,q')- \tau_r^{(s,q')}(r,q,r,s) - \tau_r^{(q,s)}(r,s,r,q')+ \tau_r^{(q,q')}(r,s,r,s)$, & $r = p = p', \; q \neq s, \; q' \neq s$, \\[1em]
						   $0$, & else, \end{tabular}\right. 
\end{split}
\end{align}
and for $(r,s) \neq (i,j)$:
\begin{align}\label{equ:asycovrsij}
\begin{split}
& \frac{\sigma_{rs,il}(p,q,p',q')}{N}\\ & \qquad= \left\{\begin{tabular}{p{5.3cm}p{5cm}} $\tau_r^{(s,l)}(p,q,p',q')$, & $r =i, \; p \notin \{i, p'\}, \; r \neq p'$, \\[1em]
                                     $-\tau_r^{(s,q')}(p,q,i,l)$, & $r =p', \; p \notin \{i, p'\}, \; r \neq i$,\\[1em]
                                     $-\tau_p^{(q,q')}(r,s,i,l)$, & $p =i, \; r \notin \{i, p'\}, \; p \neq p'$,\\[1em]
                                     $\tau_p^{(q,q')}(r,s,i,l)$, & $p =p', \; r \notin \{i, p'\}, \; p \neq i$,\\[1em]
                                     $\tau_r^{(s,l)}(p,q,r,q') - \tau_r^{(s,q')}(p,q,r,j)$, & $r=i=p', \; p \notin \{i,p'\}, \; q' \neq l$,\\[1em]
                                     $-\tau_p^{(q,j)}(r,s,p,q') + \tau_p^{(q,q')}(r,s,i,l)$, & $p=i=p', \; r \notin \{i,p'\}, \; q' \neq l$,\\[1em] 
                                     $\tau_r^{(s,l)}(r,q,p',q') - \tau_r^{(q,l)}(r,s,p',q')$, & $r=i=p, \; p' \notin \{i,p\}, \; q \neq s$,\\[1em]
                                     $-\tau_r^{(s,q')}(r,q,i,l) + \tau_r^{(q,q')}(r,s,i,l)$, & $p=r=p', \; i \notin \{r,p\}, \; q \neq s$,\\[1em]
                                     $\tau_r^{(s,l)}(p,q,p,q') + \tau_p^{(q,q')}(r,s,r,j)$, & $r=i, \; p = p', \; r\neq p', \; p \neq i$,\\[1em]
                                     $-\tau_r^{(s,q')}(p,q,p,l) - \tau_p^{(q,l)}(r,s,r,q')$, & $r=p', \; p = i, \; r \neq i, \; p \neq p'$,\\[1em]
                                     $\tau_r^{(s,l)}(r,q,r,q') - \tau_r^{(s,q')}(r,q,r,l)-\tau_r^{(q,l)}(r,s,r,q') + \tau_r^{(q,q')}(r,s,r,l)$, & $r=p=p'=i, \; s \neq q \neq q' \neq l$,\\[1em]
                                     $0$, & else, \end{tabular}\right.
\end{split}
\end{align}
where \[\tau_r^{(s,j)}(p,q,p',q')=\frac{1}{n_r}\Erw\left[\left(F_{pq}(X_{rs1})-w_{pqrs}\right)\left(F_{p'q'}(X_{rj1})-w_{p'q'rj}\right)\right].\] 

\subsection{Confidence intervals for ratios}\label{sec:tecdet_ratio}
As already shown in the main part, $\bs T_N^\text{ratio}$ is asymptotically multivariate normally distributed with expectation zero and covariance matrix $\bs S=(s_{\ell m})_{\ell,m=1}^q$ with entries
$s_{\ell m}=\frac{(\theta_\ell \bs d_\ell-\bs c_\ell)' \bs V (\theta_m \bs d_m-\bs c_m)}{\sqrt{(\theta_\ell \bs d_\ell-\bs c_\ell)' \bs V (\theta_\ell \bs d_\ell-\bs c_\ell)} \sqrt{(\theta_m \bs d_m-\bs c_m)' \bs V (\theta_m \bs d_m-\bs c_m)}}$. 
One difficulty in the construction of adequate CIs for the ratio $\theta_\ell$ is the dependency of the test statistic and the object of estimation. For the easiest case of $q=1$, one only has to deal with a single ratio 
$\theta=\frac{\bs c' \bs p}{\bs d' \bs p}$. After an application of Fieller's (1954)\nocite{fieller1954some} theorem, a two-sided CI follows by solving the inequality
\begin{equation}\label{equ:CIrat}
\frac{|\sqrt{N}\left(\bs\theta \bs d- \bs c \right)'\left(\widehat{\bs p}-\bs p\right)|}{\bs \theta^2 \bs d' \widehat{\bs V}_N \bs d- 2 \bs\theta \bs c' \widehat{\bs V}_N \bs d + \bs c' \widehat{\bs V}_N \bs c} \leq c(\alpha),
\end{equation}
where $c(\alpha)$ is an adequate $(1-\alpha)$ equicoordinate quantile, which depends on $\tau_\ell$ through the covariance matrix $\bs S$. A way to calculate this corresponding quantile is given at the end of this section.
Another representation of inequality \eqref{equ:CIrat} is given by 
\begin{equation}\label{equ:CIratalt}
 A\bs\theta^2 + B\bs\theta + C \leq 0,
\end{equation}
where \begin{eqnarray*}
       A &=& \left(\sqrt{N}\bs d'\widehat{\bs p}\right)^2-\left(c(\alpha)\right)^2 \bs d'\widehat{\bs V}_N \bs d, \\
       B &=& -2\left[\left(\sqrt{N}\bs c'\widehat{\bs p}\right)\left(\sqrt{N}\bs d'\widehat{\bs p}\right)-\left(c(\alpha)\right)^2 \bs c'\widehat{\bs V}_N \bs d\right]\; \text{ and} \; \\
       C &=& \left(\sqrt{N}\bs c'\widehat{\bs p}\right)^2-\left(c(\alpha)\right)^2 \bs c'\widehat{\bs V}_N \bs c.
      \end{eqnarray*}
In case of $q=1$, three possible solutions of inequality \eqref{equ:CIratalt} exist. The first one results in the case that all values lie outside the finite interval defined by the two roots of inequalities, the second one results in the entire $\theta$-axis.
The third solution results in a finite interval which is the most desirable solution in this case. This solution corresponds to $A >0$ and thus, $B^2-4AC > 0$ since the first two solutions only occur with small probability if $\bs d' \bs p$ is significantly from 
zero. Dilba et al.~(2006)\nocite{dilba2006simultaneous} gave another representation of $A>0$, namely $\frac{(c(\alpha))^2 \bs d'\widehat{\bs V}_N \bs d}{\left(\sqrt{N}\bs d'\widehat{\bs p}\right)^2}>1$. For the general ratio problem 
with arbitrary $q \in \mathbb{N}$ and in order to guarantee 
$A >0$ for each component $\ell=1,\ldots,q$, say $A_\ell>0$, with high probability the following inequalities must be fulfilled
\begin{equation}
 0 < \frac{y\sqrt{\bs d'\widehat{\bs V}_N \bs d}}{\sqrt{N}\bs d'\widehat{\bs p}} \ll 1  \quad \text{or} \quad \frac{\sqrt{N}\bs d'\widehat{\bs p}}{y\sqrt{\bs d'\widehat{\bs V}_N \bs d}} \gg 1
\end{equation}
for some relevant point $y$ (see Dilba et al., 2006)\nocite{dilba2006simultaneous}. 

To evaluate the right $(1-\alpha)$ equicoordinate quantile for deriving the corresponding SCI, Dilba et al.~(2006)\nocite{dilba2006simultaneous} introduce three different approaches. Here, we only focus on a resampling approach 
based on the wild bootstrap introduced in Sections \ref{sec:boot} and \ref{sec:MCTPboot}. A wild bootstrap based CI can be determined by the following algorithm:
\begin{enumerate}
 \item For $b=1,\ldots,B$:
  \[\text{Calculate the wild bootstrap version $T_\ell^{\epsilon,b}$ of the corresponding test statistic $T_\ell$.}\]
 \item Compute the $(1-\alpha)$ quantile of the values: \[T_\text{max}^{\epsilon,b}=\max\{|T_1^{\epsilon,b}|, \ldots, |T_q^{\epsilon,b}|\} \; \text{for}\;  b=1,\ldots,B.\]
\end{enumerate}
The result of the computation in step 2 leads to the $(1-\alpha)$ equicoordinate wild bootstrap based quantile $c^\epsilon(\alpha)$. Finally, combining all the latter results a wild bootstrap based SCI for the ratio $\theta_\ell$ can be calculated by
\[\left[\frac{-B_\ell - \sqrt{B_\ell^2-4A_\ell C_\ell}}{2A_\ell}, \frac{-B_\ell + \sqrt{B_\ell^2-4A_\ell C_\ell}}{2A_\ell}\right],\]
where \begin{eqnarray*}
       A_\ell &=& \left(\sqrt{N}\bs d'_\ell\widehat{\bs p}\right)^2-\left(c^\epsilon(\alpha)\right)^2 \bs d'_\ell\widehat{\bs V}_N \bs d_\ell, \\
       B_\ell &=& -2\left[\left(\sqrt{N}\bs c'_\ell\widehat{\bs p}\right)\left(\sqrt{N}\bs d'_\ell\widehat{\bs p}\right)-\left(c^\epsilon(\alpha)\right)^2 \bs c'_\ell\widehat{\bs V}_N \bs d_\ell\right]\; \text{ and} \; \\
       C_\ell &=& \left(\sqrt{N}\bs c'_\ell\widehat{\bs p}\right)^2-\left(c^\epsilon(\alpha)\right)^2 \bs c'_\ell\widehat{\bs V}_N \bs c_\ell.
      \end{eqnarray*}

\section{The Proofs}\label{sec:proof}

\subsection*{Proof of Theorem~\ref{asy1:rm}}\label{app_rm}

First, note that given the data, only the Rademacher variables $\epsilon_{ik}$ are random and all other quantities are deterministic. Especially the $D_{ijk}(r,s)$ are deterministic sequences. Second, note that $\sqrt{N}\widehat{\bs w}^\epsilon$ 
can be represented as a continuous function of the pooled $(i=1,\ldots,a)$ random vectors
\[\sqrt{N}\frac{1}{n_i}\sum\limits_{k=1}^{n_i}\epsilon_{ik}\left(D_{ijk}(r,s)\right)_{r,s=1}^{a,d},\]
which are just sums of row-wise independent random vectors.

The result follows from an application of a multivariate Lindeberg-Feller theorem. Therefore, see Theorem A.1 in Beyersmann et al.~(2013)\nocite{beyersmann2013weak} and Theorem 4.1 in Pauly (2011)\nocite{pauly2011weighted}. 
First, we have to show that the quantity $\frac{1}{\sqrt{N}}\epsilon_{ik}(D_{ijk}(r,s))_{r,s=1}^{a,d}$ fulfilles the multivariate Lindeberg condition given the data. Defining $ c\sqrt{n_i}=:s_{n_i}$, it can be easily shown that
\[\frac{1}{s_{n_i}^2}\sum_{k=1}^{n_i}\Erw\left((\epsilon_{ik}D_{ijk}(r,s))^2\mathbbm{1}\{|\epsilon_{ik}D_{ijk}(r,s)|>\epsilon s_{n_i}^2\}\right)\stackrel{n_i\to\infty}{\to}0,\] by an application of the 
dominated convergence theorem since $|\epsilon_{ik}D_{ijk}(r,s)|\leq 1$. 
Furthermore, conditioned on the 
data the expectation of $\sqrt{N}\hat{w}_{pqrs}^\epsilon$ is zero. Thus, it remains to show that the conditional covariance between $\sqrt{N}\hat{w}_{pqrs}^\epsilon$ and $\sqrt{N}\hat{w}_{p'q'ij}^\epsilon$ converges in probability to the particular 
element of $\bm{\Sigma}$ for each $p,q,r,s,p',q',i,j$. First, the 
calculations regarding \[\sigma_{rs}(p, q, p', q') = N \cdot \Cov(\sqrt{N}\hat{w}_{pqrs}^\epsilon , \sqrt{N}\hat{w}_{p'q'rs}^\epsilon|\bs X)\] are given. For $r \neq p, \; r \neq p', \; p \neq p'$ it follows due to independence of the Rademacher variables:
\begin{align*}
 &\Cov\left(\sqrt{N}\hat{w}^\epsilon_{pqrs}, \sqrt{N}\hat{w}^\epsilon_{p'q'rs} \big\vert \bs X \right) \\ &\qquad = N \cdot \Cov\left(\frac{1}{n_r}\sum_{k=1}^{n_r} \epsilon_{rk}D_{rsk}(p,q)-\frac{1}{n_p}\sum_{k=1}^{n_p} \epsilon_{pk}D_{pqk}(r,s), \right.\\
													&\qquad \; \; \; \left. \vphantom{\frac{1}{n_r}\sum_{k=1}^{n_r}}
													\frac{1}{n_r}\sum_{k=1}^{n_r} \epsilon_{rk}D_{rsk}(p',q')-\frac{1}{n_{p'}}\sum_{k=1}^{n_{p'}} \epsilon_{p'k}D_{p'q'k}(r,s)\bigg\vert\bs X\right)\\
										  &\qquad=  N \cdot \Cov\left(\frac{1}{n_r}\sum_{k=1}^{n_r}\epsilon_{rk}D_{rsk}(p,q), \; \frac{1}{n_r}\sum_{k=1}^{n_r}\epsilon_{rk}D_{rsk}(p',q')\bigg\vert\bs X\right) \\
										  &\qquad= \frac{N}{n_r^2}\sum_{k=1}^{n_r}D_{rsk}(p,q)D_{rsk}(p',q')\\
										  &\qquad\stackrel{p}{\rightarrow} \tau_r^{(s,s)}(p,q,p',q').
\end{align*}
For $r \neq p, \; r \neq p', \; p = p'$ we calculate:
\begin{align*}
& \Cov\left(\sqrt{N}\hat{w}^\epsilon_{pqrs}, \sqrt{N}\hat{w}^\epsilon_{pq'rs} \big\vert \bs X \right)\\ &\qquad= N \cdot \Cov\left(\frac{1}{n_r}\sum_{k=1}^{n_r} \epsilon_{rk}D_{rsk}(p,q)-\frac{1}{n_p}\sum_{k=1}^{n_p} \epsilon_{pk}D_{pqk}(r,s), \right. \\
													&\qquad \; \; \; \left. \vphantom{\frac{1}{n_r}\sum_{k=1}^{n_r}}  
													\frac{1}{n_r}\sum_{k=1}^{n_r} \epsilon_{rk}D_{rsk}(p,q')-\frac{1}{n_p}\sum_{k=1}^{n_p} \epsilon_{pk}D_{pq'k}(r,s)\bigg\vert\bs X\right)\\
										  &\qquad=  N \cdot \Cov\left(\frac{1}{n_r}\sum_{k=1}^{n_r}\epsilon_{rk}D_{rsk}(p,q), \; \frac{1}{n_r}\sum_{k=1}^{n_r}\epsilon_{rk}D_{rsk}(p,q')\bigg\vert\bs X\right) \\
										  &\qquad \; \; \; 
										  +N \cdot \Cov\left(-\frac{1}{n_p}\sum_{k=1}^{n_p}\epsilon_{pk}D_{pqk}(r,s), \; -\frac{1}{n_p}\sum_{k=1}^{n_p}\epsilon_{pk}D_{pq'k}(r,s)\bigg\vert\bs X\right)\\
										  &\qquad= \frac{N}{n_r^2}\sum_{k=1}^{n_r}D_{rsk}(p,q)D_{rsk}(p,q') + \frac{N}{n_p^2}\sum_{k=1}^{n_p}D_{pqk}(r,s)D_{pq'k}(r,s)\\
										  &\qquad\stackrel{p}{\rightarrow} \tau_r^{(s,s)}(p,q,p,q') + \tau_{p}^{(q,q')}(r,s,r,s).
\end{align*}

Similar, for $r = p, \; p \neq p', \; q \neq s$, we obtain:
\begin{align*}
& \Cov\left(\sqrt{N}\hat{w}^\epsilon_{rqrs}, \sqrt{N}\hat{w}^\epsilon_{p'q'rs} \big\vert \bs X \right)\\ &\qquad= N \cdot \Cov\left(\frac{1}{n_r}\sum_{k=1}^{n_r} \epsilon_{rk}D_{rsk}(r,q)-\frac{1}{n_r}\sum_{k=1}^{n_r} \epsilon_{rk}D_{rqk}(r,s), \right. \\
													&\qquad \; \; \; \left. \vphantom{\frac{1}{n_r}\sum_{k=1}^{n_r}}  
													\frac{1}{n_r}\sum_{k=1}^{n_r} \epsilon_{rk}D_{rsk}(p',q')-\frac{1}{n_{p'}}\sum_{k=1}^{n_{p'}} \epsilon_{p'k}D_{p'q'k}(r,s)\bigg\vert\bs X\right)\\
										  &\qquad= N \cdot Cov\left(\frac{1}{n_r}\sum_{k=1}^{n_r}\epsilon_{rk}D_{rsk}(r,q)-\frac{1}{n_r}\sum_{k=1}^{n_r}\epsilon_{rk}D_{rqk}(r,s),% \right. \\
										      %&\qquad\;\;\; \left. \vphantom{\frac{1}{n_r}\sum_{k=1}^{n_r}} 
										      \frac{1}{n_r}\sum_{k=1}^{n_r}\epsilon_{rk}D_{rsk}(p',q')\bigg\vert\bs X\right) \\
										  &\qquad=  N \cdot \Cov\left(\frac{1}{n_r}\sum_{k=1}^{n_r}\epsilon_{rk}D_{rsk}(r,q), \; \frac{1}{n_r}\sum_{k=1}^{n_r}\epsilon_{rk}D_{rsk}(p',q')\bigg\vert\bs X\right) \\
										  &\qquad \; \; \; 
										  +N \cdot \Cov\left(-\frac{1}{n_r}\sum_{k=1}^{n_r}\epsilon_{rk}D_{rqk}(r,s), \; \frac{1}{n_r}\sum_{k=1}^{n_r}\epsilon_{rk}D_{rsk}(p',q')\bigg\vert\bs X\right)\\
										  &\qquad= \frac{N}{n_r^2}\sum_{k=1}^{n_r}D_{rsk}(r,q)D_{rsk}(p',q') - \frac{N}{n_r^2}\sum_{k=1}^{n_r}D_{rqk}(r,s)D_{rsk}(p',q')\\
										  &\qquad\stackrel{p}{\rightarrow} \tau_r^{(s,s)}(r,q,p',q') - \tau_{r}^{(q,s)}(r,s,p',q'),
\end{align*}
for $r = p', \; p \neq p', \; q' \neq s$:
\begin{align*}
& \Cov\left(\sqrt{N}\hat{w}^\epsilon_{pqrs}, \sqrt{N}\hat{w}^\epsilon_{rq'rs} \big\vert \bs X \right)\\ &\qquad= N \cdot \Cov\left(\frac{1}{n_r}\sum_{k=1}^{n_r} \epsilon_{rk}D_{rsk}(p,q)-\frac{1}{n_p}\sum_{k=1}^{n_p} \epsilon_{pk}D_{pqk}(r,s), \right. \\
													&\qquad \; \; \; \left. \vphantom{\frac{1}{n_r}\sum_{k=1}^{n_r}} 
													\frac{1}{n_r}\sum_{k=1}^{n_r} \epsilon_{rk}D_{rsk}(r,q')-\frac{1}{n_r}\sum_{k=1}^{n_r} \epsilon_{rk}D_{rq'k}(r,s)\bigg\vert\bs X\right)\\
										  &\qquad= N \cdot Cov\left(\frac{1}{n_r}\sum_{k=1}^{n_r}\epsilon_{rk}D_{rsk}(p,q), %\right. \\
										      %&\qquad\;\;\; \left. \vphantom{\frac{1}{n_r}\sum_{k=1}^{n_r}} 
										      \frac{1}{n_r}\sum_{k=1}^{n_r}\epsilon_{rk}D_{rsk}(r,q')-\frac{1}{n_r}\sum_{k=1}^{n_r}\epsilon_{rk}D_{rq'k}(r,s)\bigg\vert\bs X\right) \\
										  &\qquad=  N \cdot \Cov\left(\frac{1}{n_r}\sum_{k=1}^{n_r}\epsilon_{rk}D_{rsk}(p,q), \; \frac{1}{n_r}\sum_{k=1}^{n_r}\epsilon_{rk}D_{rsk}(r,q')\bigg\vert\bs X\right) \\
										  &\qquad \; \; \; 
										  +N \cdot \Cov\left(\frac{1}{n_r}\sum_{k=1}^{n_r}\epsilon_{rk}D_{rsk}(p,q), \; -\frac{1}{n_r}\sum_{k=1}^{n_r}\epsilon_{rk}D_{rq'k}(r,s)\bigg\vert\bs X\right)\\
										  &\qquad= \frac{N}{n_r^2}\sum_{k=1}^{n_r}D_{rsk}(p,q)D_{rsk}(r,q') - \frac{N}{n_r^2}\sum_{k=1}^{n_r}D_{rsk}(p,q)D_{rq'k}(r,s)\\
										  &\qquad\stackrel{p}{\rightarrow} \tau_r^{(s,s)}(p,q,r,q') - \tau_{r}^{(s,q')}(p,q,r,s),
\end{align*}
and finally, for $r = p, \; r =p', \; q \neq s,  \; q' \neq s$:
\begin{align*}
& \Cov\left(\sqrt{N}\hat{w}^\epsilon_{rqrs}, \sqrt{N}\hat{w}^\epsilon_{rq'rs} \big\vert \bs X \right)\\ &\qquad= N \cdot \Cov\left(\frac{1}{n_r}\sum_{k=1}^{n_r} \epsilon_{rk}D_{rsk}(r,q)-\frac{1}{n_r}\sum_{k=1}^{n_r} \epsilon_{rk}D_{rqk}(r,s), \right. \\
													&\qquad \; \; \; \left. \vphantom{\frac{1}{n_r}\sum_{k=1}^{n_r}}  
													\frac{1}{n_r}\sum_{k=1}^{n_r} \epsilon_{rk}D_{rsk}(r,q')-\frac{1}{n_r}\sum_{k=1}^{n_r} \epsilon_{rk}D_{rq'k}(r,s)\bigg\vert\bs X\right)\\
										  &\qquad=  N \cdot \Cov\left(\frac{1}{n_r}\sum_{k=1}^{n_r}\epsilon_{rk}D_{rsk}(r,q), \; \frac{1}{n_r}\sum_{k=1}^{n_r}\epsilon_{rk}D_{rsk}(r,q')\bigg\vert\bs X\right) \\
										  &\qquad \; \; \; 
										  +N \cdot \Cov\left(-\frac{1}{n_r}\sum_{k=1}^{n_r}\epsilon_{rk}D_{rqk}(r,s), \; \frac{1}{n_r}\sum_{k=1}^{n_r}\epsilon_{rk}D_{rsk}(r,q')\bigg\vert\bs X\right)\\
										  &\qquad \; \; \; +N \cdot \Cov\left(\frac{1}{n_r}\sum_{k=1}^{n_r}\epsilon_{rk}D_{rsk}(r,q), \; -\frac{1}{n_r}\sum_{k=1}^{n_r}\epsilon_{rk}D_{rq'k}(r,s)\bigg\vert\bs X\right)\\
										  &\qquad \; \; \; 
										  +N \cdot \Cov\left(-\frac{1}{n_r}\sum_{k=1}^{n_r}\epsilon_{rk}D_{rqk}(r,s), \; -\frac{1}{n_r}\sum_{k=1}^{n_r}\epsilon_{rk}D_{rq'k}(r,s)\bigg\vert\bs X\right)\\
										  &\qquad= \frac{N}{n_r^2}\sum_{k=1}^{n_r}D_{rsk}(r,q)D_{rsk}(r,q') - \frac{N}{n_r^2}\sum_{k=1}^{n_r}D_{rqk}(r,s)D_{rsk}(r,q')\\
										  &\qquad \; \; \; 
										  -\frac{N}{n_r^2}\sum_{k=1}^{n_r}D_{rsk}(r,q)D_{rq'k}(r,s)+ \frac{N}{n_r^2}\sum_{k=1}^{n_r}D_{rqk}(r,s)D_{rq'k}(r,s) \\
										  &\qquad\stackrel{p}{\rightarrow} \tau_r^{(s,s)}(r,q,r,q') - \tau_r^{(q,s)}(r,s,r,q') - \tau_r^{(s,q')}(r,q,r,s) + \tau_{r}^{(q,q')}(r,s,r,s).
\end{align*}
In all other cases the conditional variance is 0. 

Concerning the case of the covariance, where $(r,s) \neq (i,j)$, one can proceed in a similar way. Here, only the case $r=i, \;p \neq i, \;p \neq p', \; r \neq p'$ is shown as an example, since the 
calculation is straightforward as in case of the variance given previously:
\begin{align*}
& \Cov\left(\sqrt{N}\hat{w}^\epsilon_{pqrs}, \sqrt{N}\hat{w}^\epsilon_{p'q'rj} \big\vert \bs X \right)\\ &\qquad= N \cdot \Cov\left(\frac{1}{n_r}\sum_{k=1}^{n_r} \epsilon_{rk}D_{rsk}(p,q)-\frac{1}{n_p}\sum_{k=1}^{n_p} \epsilon_{pk}D_{pqk}(r,s), \right. \\
													&\qquad \; \; \; \left. \vphantom{\frac{1}{n_r}\sum_{k=1}^{n_r}}
													\frac{1}{n_r}\sum_{k=1}^{n_r} \epsilon_{rk}D_{rjk}(p',q')-\frac{1}{n_p'}\sum_{k=1}^{n_p'} \epsilon_{p'k}D_{p'q'k}(r,j)\bigg\vert\bs X\right)\\
										  &\qquad=  N \cdot \Cov\left(\frac{1}{n_r}\sum_{k=1}^{n_r}\epsilon_{rk}D_{rsk}(p,q), \; \frac{1}{n_r}\sum_{k=1}^{n_r}\epsilon_{rk}D_{rjk}(p',q')\bigg\vert\bs X\right) \\
										  &\qquad= \frac{N}{n_r^2}\sum_{k=1}^{n_r}D_{rsk}(p,q)D_{rjk}(p',q')\\
										  &\qquad\stackrel{p}{\rightarrow} \tau_r^{(s,j)}(p,q,p',q').
\end{align*}
All other cases can be treated similarly and finally, it follows that the covariance of $\sqrt{N}\widehat{\bs w}^\epsilon$ conditioned on the data converges in probability to the covariance matrix $ \bm{\Sigma}$. All in all, the distribution of 
$\sqrt{N}\widehat{\bs w}^\epsilon$ weakly converges to a multivariate normal distribution with mean zero and covariance matrix $ \bm{\Sigma}$ in probability. 

And again by the following equation $\widehat{\bs p}=\bs E_{ad} \cdot \widehat{\bs w}$, it follows that $\sqrt{N}\widehat{\bs p}^\epsilon$ is asymptotically multivariate normally distributed with expectation $\bs 0$ and asymptotic covariance
matrix $\bs V_N=\bs E_{ad}\bm{\Sigma} \bs E_{ad}$. 
\hfill$\square$

\subsection*{Proof of Corollary \ref{cor:ATS}}
Using the asymptotic distribution under $H_0^p$ of $F_N(\bs M)$, which is given by 
\[F(\bs M)=\sum_{i=1}^a \sum_{j=1}^d \frac{\lambda_{ij}(\bs M \bs V_N)}{\tr(\bs M \bs V_N)} C_{ij}^2,\] where $C_{ij}$ are independent standard normal random variables and $\lambda_{ij}(\bs{MV})$ denote the corresponding eigenvalues of $\bs{MV}$,
%The authors propose to approximate the distribution of $F_N(\bs M)$ by 
the results follow immediately from Theorem \ref{asy1:rm} by applying the continuous mapping theorem and Slutsky's theorem: \[N\widehat{\bs p}' \bs M \widehat{\bs p}\stackrel{d}{\to} A \sim \sum_{i=1}^a\sum_{j=1}^d\lambda_{ij}(\bs{MV}) \xi_{ij},\]
where $\xi_{ij}\stackrel{i.i.d}{\sim}\chi^2_1$ and the multivariate normal distribution of $\sqrt{N}\widehat{\bs p}^\epsilon$.
\hfill$\square$

\subsection*{Proof of Corollary \ref{cor:distr}}
Using definition of a contrast matrix $\bs C=(\bs c_1 ,\ldots ,\bs c_q)'$, the proof directly follows by the asymptotic normality of $\sqrt{N}\bs C \widehat{\bs p}^\epsilon$ and an application of Slutsky's theorem.
\hfill$\square$

\section{Power Simulations}
\begin{figure}[ht]
\centering
   \includegraphics[width=.9\linewidth]{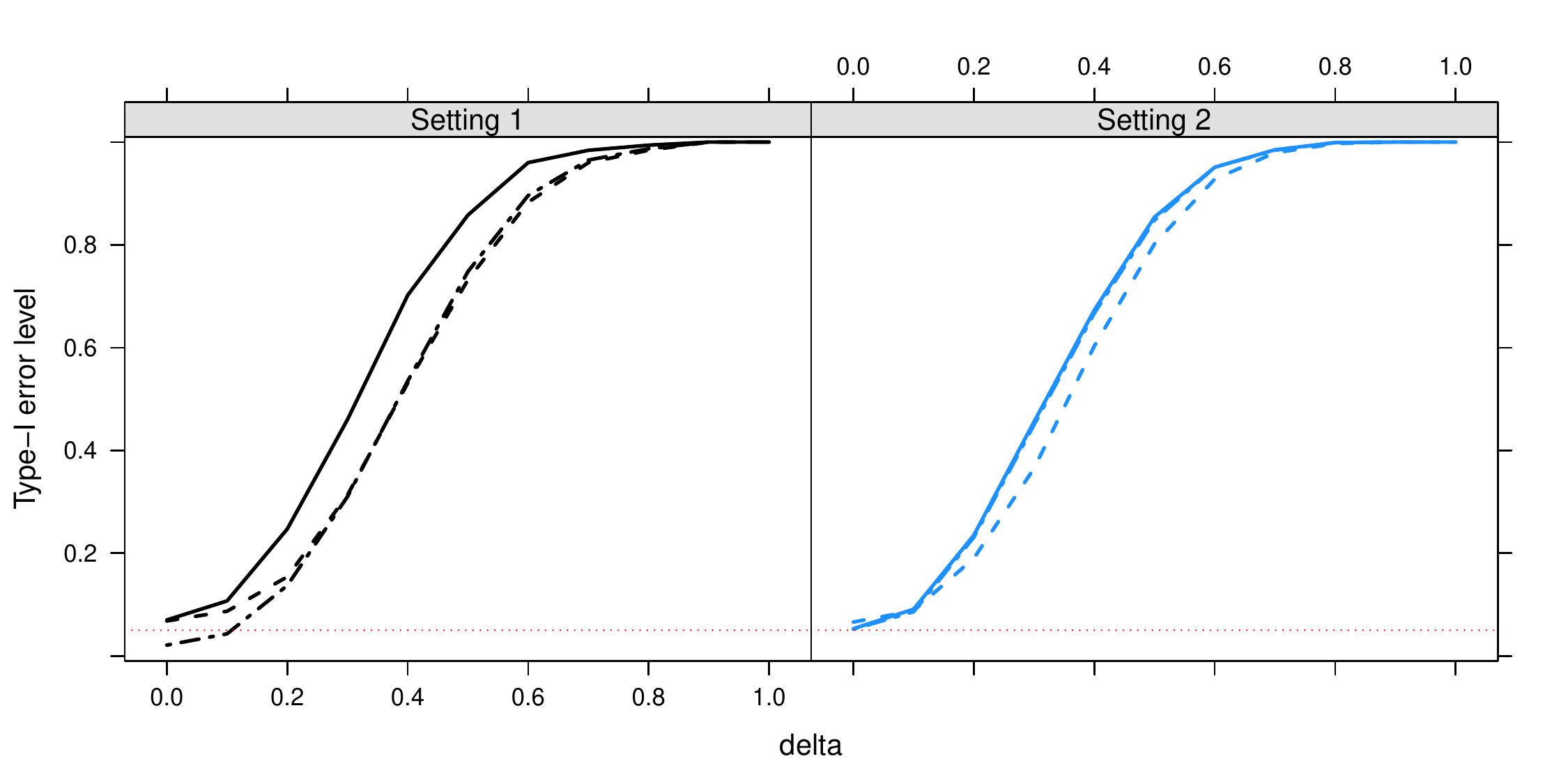}
  \caption{\label{fig:powerA}Power simulation results for Setting 1 (left) and Setting 2 (right) for factor $A$ with a compound symmetry covariance structure for three different tests, namely MCTP (dotted-dashed), bootMCTP (solid) and ATS (dashed).}
  \end{figure}

\begin{figure}[ht]
\centering
  \includegraphics[width=.9\linewidth]{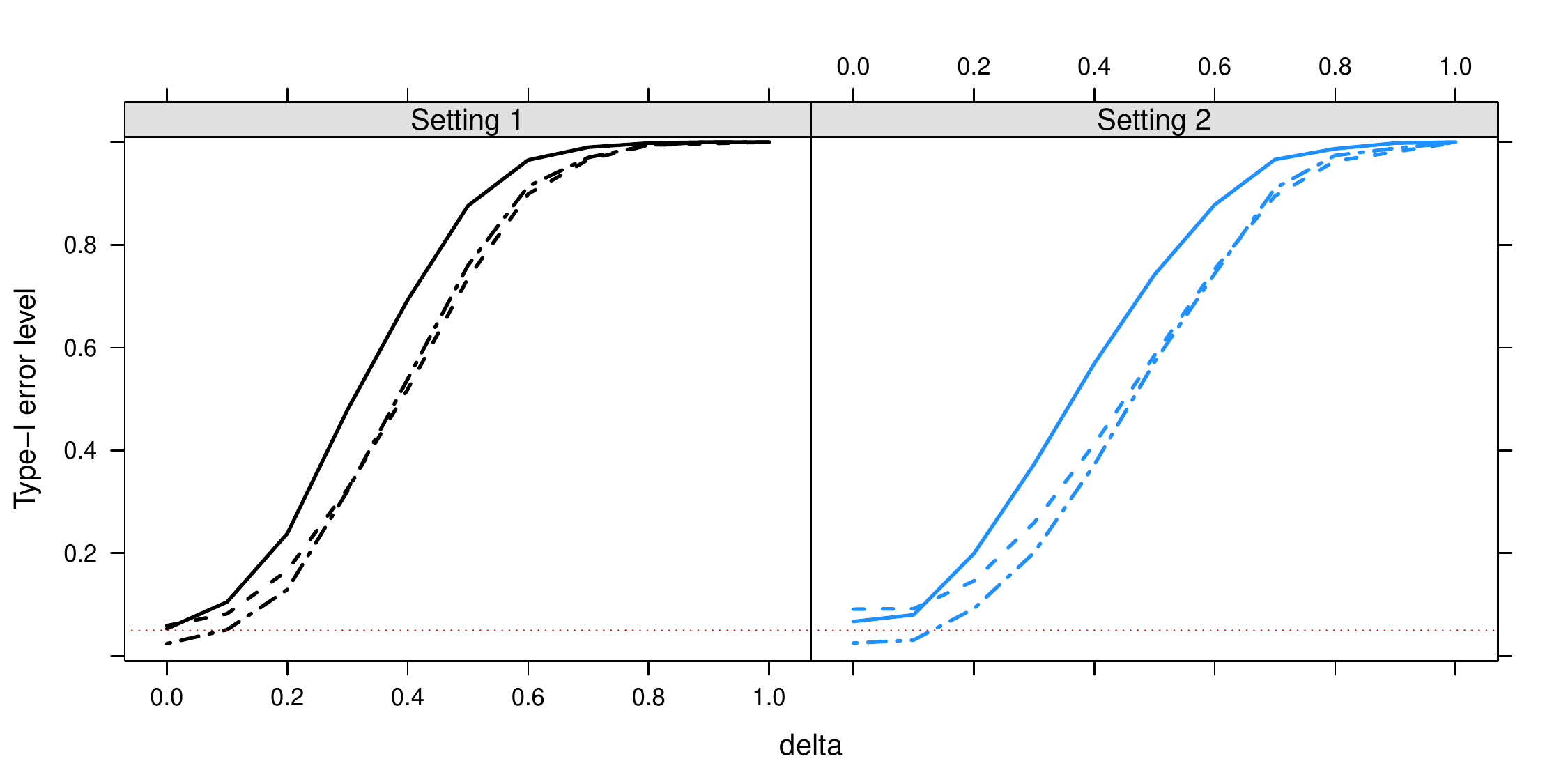}
  \caption{\label{fig:powerD}Power simulation results for Setting 1 (left) and Setting 2 (right) for factor $D$ with a compound symmetry covariance structure for three different tests, namely MCTP (dotted-dashed), bootMCTP (solid) and ATS (dashed).}
\end{figure}

In this section, a small power simulation is presented. The simulations are restricted to one covariance structure (compound symmetry) and to a centring matrix $\bs P$ as a contrast matrix. Again, Setting 1 and Setting 2 as in the main part are simulated
and a balanced, homoscedastic design
with $n_{ij}=20$ individuals per group is examined. To conduct a power simulation, an effect $\delta \in \{0, 0.1, 0.2, \ldots, 1\}$ is added to the different factors. Figure~\ref{fig:powerA} summarizes the results for an effect in factor $A$, whereas
Figure~\ref{fig:powerD} shows the results for an effect in factor $D$. In case of Setting 1 and an effect  in factor $A$, $\bm{\delta}=\left(0,0,0,0,0,0,\delta, \delta, \delta\right)'$ and regarding an effect in factor $B$, the vector is given as follows
$\bm{\delta}=\left(0,0,\delta,0,0,\delta,0,0,\delta\right)'$. For Settting 2, the $\bm{\delta}$-vector is constructed in a similar way.

For the power behavior in the first setting and main effect $A$, it can be readily seen from Figure~\ref{fig:powerA} that the bootMCTP has the highest power in case of three repeated measures when compared with the ATS and the standard MCTP. The power 
functions of the ATS and the standard MCTP are hardly distinguishable. Regarding Setting 2, the power behaviors of both MCTPs are slightly better compared to the ATS.

Regarding the results summarized in Figure~\ref{fig:powerD} and Setting 1, the same behavior as given above can be observed. Again, the wild bootstrap version of the MCTP shows the highest power. For small vales of delta the ATS is slightly better 
than the standard MCTP, but these two power curves become indistinguishable for higher values of delta. The results of Setting 2 are comparable to the results of Setting 1. 
%   \begin{figure}[ht]
% \centering
%   \includegraphics[width=.9\linewidth]{power_CS_I.pdf}
%   \caption{\label{fig:powerI}Power simulation results for Setting 1 (left) and Setting 2 (right) for factor $I$ with a compound symmetry covariance structure for three different tests, namely MCTP (dotted-dashed), bootMCTP (solid) and ATS (dashed).}
% \end{figure}
%   \textcolor{red}{Beschreibung für Figure~\ref{fig:powerI} einfügen!}
% Regarding Setting 1, the results presented in Figure~\ref{fig:powerD} exhibit the same behavior to that given above. 
% In case of Setting 2, the shape of the power curves is comparable to Setting 1 with the difference that the power in the design with four repeated measures is less than in the design with three repeated measures.

\end{document}